\documentclass[final,hidelinks]{siamart251216}

\usepackage[T1]{fontenc}
\usepackage{amsmath,amssymb,bm}
\usepackage{booktabs}
\usepackage{adjustbox}
\usepackage{graphicx}
\usepackage{siunitx}
\usepackage{arydshln}
\usepackage{placeins}
\usepackage{float}
\usepackage{xcolor}

\graphicspath{{paper_figures_ver7/}{Figures/}}

\newcommand{\ii}{\mathrm{i}}
\newcommand{\uh}{u_h}
\newcommand{\vh}{v_h}
\newcommand{\uref}{u^{\rm ref}}
\newcommand{\ppw}{\mathrm{ppw}}
\newcommand{\npml}{n_{\rm pml}}

\newcommand{\revision}[1]{#1}

\headers{A Massively Parallel Three-Grid Preconditioner}%
{S. Fu, Y. Wang, and Z. Zhao}
\begin{document}

\title{A Massively Parallel Three-Grid Preconditioner for the High-Frequency
Helmholtz Equation}
\author{Shubin Fu\thanks{Corresponding author. School of Mathematical Sciences,
Eastern Institute of Technology, Ningbo, Zhejiang 315200, P. R. China
(\email{sfu@eitech.edu.cn}).}
\and Yitong Wang\thanks{School of Mathematical Sciences, Shanghai Jiao Tong
University, Shanghai 200240, P. R. China; School of Mathematical Sciences, Eastern
Institute of Technology, Ningbo, Zhejiang 315200, P. R. China.}
\and Zixiao Zhao\footnotemark[2]}
\maketitle

\begin{abstract}
Accurate simulation of three-dimensional time-harmonic wave propagation over
many wavelengths requires control of phase error and efficient solution of
large indefinite systems.  We develop a three-grid solver based on the compact
27-point interpolated optimized finite-difference (IOFD) discretization.  Its
wavenumber-dependent stencil supports a fine-grid
resolution of six points per shortest wavelength and an unshifted physical
correction on the \(2h\) grid at only three points per shortest wavelength.  The
method retains unshifted IOFD operators on the \(h\) and \(2h\) grids, while a
complex-shifted \(2h\)--\(4h\) auxiliary cycle preconditions a
factorization-free iterative approximation of the coarse inverse.  Restricting
the shift to this auxiliary cycle preserves the propagative character of the
coarse correction.  Comparison with the outgoing Green function confirms phase
and relative-amplitude accuracy on a sequence of meshes up to \revision{\(13376^3\)},
\revision{exceeding two trillion unknowns,} with
the largest problem spanning approximately \revision{2227} wavelengths per coordinate.  The
same fixed solver configuration retains
robust convergence across smooth, discontinuous, high-contrast, and geophysical
velocity models and exhibits scalable parallel performance.  In particular, a
\revision{\(2560^3\)} problem spanning approximately \revision{425} wavelengths in each coordinate direction is
solved in \revision{21.7 seconds} on just 64 NVIDIA A100 GPUs.
\end{abstract}



\section{Introduction}

Large-scale frequency-domain wave simulation is constrained by two coupled
difficulties: phase accuracy and linear-solver complexity.  In seismic imaging
and inverse scattering, the Helmholtz equation is solved repeatedly for
many sources and frequencies in heterogeneous domains spanning hundreds of
wavelengths \cite{VirieuxOperto2009,TournierEtAl2022}.  At fixed resolution, a
three-dimensional domain containing \(K\) wavelengths in each coordinate
direction has \(\mathcal{O}(K^3)\) unknowns, while truncation by a perfectly
matched layer (PML) produces a complex non-Hermitian indefinite operator.
Numerical phase error also accumulates with propagation distance and produces
the pollution effect \cite{BabuskaSauter1997}.  An effective high-frequency
method must therefore control dispersion without making the resulting linear
systems prohibitively expensive to solve.

Compact optimized finite differences address the dispersion side of this
problem by matching the discrete and continuous dispersion relations
\cite{JoShinSuh1996,AghamiryEtAl2022,CocquetGander2024}.  In particular, the
interpolated optimized finite-difference (IOFD) scheme of Stolk uses a compact
27-point stencil whose coefficients vary smoothly with the local
nondimensional wavenumber \(kh\) \cite{Stolk2016IOFD}.  Its directional
phase-slowness optimization accurately represents wave propagation at
approximately five to six points per wavelength and remains effective near
three points per wavelength.  IOFD can therefore reduce the target-grid problem
size while retaining a compact stencil, and it makes a geometrically coarsened
physical operator plausible even when the \(2h\) grid lies close to its sampling
limit.

Once the discretization has been fixed, the remaining task is to solve a large
indefinite system, often for many right-hand sides at each frequency.  Sparse
direct solvers are attractive in this setting because one factorization can be
reused across all sources, leaving comparatively inexpensive triangular solves
\cite{OpertoEtAl2007}.  In three dimensions, however, fill-in causes both
factorization time and memory to grow much faster than the number of unknowns,
so the setup becomes prohibitive as the frequency and domain size increase.

Iterative Helmholtz solvers typically employ Krylov methods whose effectiveness
depends on preconditioners that approximate long-range wave propagation through
reduced or decomposed problems.  Multilevel Schwarz methods use
optimized transmission conditions together with Dirichlet-to-Neumann or
spectral coarse spaces
\cite{GanderMagoulesNataf2002,GrahamSpenceVainikko2017,ConenEtAl2014,
BootlandEtAl2021,TournierEtAl2022}, whereas sweeping preconditioners approximate
block factorizations or sequences of Dirichlet-to-Neumann maps
\cite{EngquistYing2011,Stolk2013DD,GanderZhang2019,PoulsonEtAl2013}.
Travel-time factorizations, spectral deflation, and learned corrections provide
other representations of global wave error
\cite{TreisterHaber2019,AzulayTreister2023,LererBenYairTreister2024,
ChenDwarkaVuik2025}.  WaveHoltz takes a different route by filtering a
periodically forced wave equation, with a recent overset-grid realization
attaining linear complexity at fixed frequency
\cite{AppeloGarciaRunborg2020,AppeloEtAl2026}.  These approaches demonstrate the
value of incorporating global wave information.  Their realizations may,
however, require local or panel factorizations, global coarse solves, ordered
traversals, or an auxiliary time-domain evolution; low-rank Green-function
representations also become more expensive with frequency
\cite{EngquistZhao2018}.

Geometric multigrid is particularly attractive for compact finite differences
because its level operators, relaxation, and transfer can all be applied through
local stencil and vector operations.  Its use for the Helmholtz equation is
difficult, however, because aggressive coarsening changes the character of the
discrete waves.  Classical multigrid deteriorates when relaxation and coarse
correction fail to represent oscillatory error \cite{GanderZhang2019}.  The
complex-shifted Laplacian and related two-grid methods stabilize the hierarchy
through artificial absorption
\cite{ErlanggaOosterleeVuik2006,VanGijzenErlanggaVuik2007,CalandraEtAl2013}, but
a shifted physical coarse operator no longer reproduces the propagation carried
by the fine grid.  Wave-ray multigrid addresses coarse representation by
transferring slowly varying ray amplitudes \cite{BrandtLivshits1997}, while
dispersion-matched operators preserve the fine-grid phase relation under
aggressive coarsening
\cite{StolkAhmedBhowmik2014,Stolk2016IOFD,Stolk2025TwoGrid}.  Local Fourier
analysis (LFA) provides a quantitative framework for comparing such
rediscretizations and their relaxation schemes \cite{YovelTreister2024}.

Taken together, these multigrid developments separate two roles of the coarse
level.  The physical coarse operator must represent wave propagation, whereas
its inverse must admit an efficient approximation.  The regime considered here
requires both properties at once: the \(2h\) grid has only three points per
shortest wavelength, and the entire hierarchy must remain factorization-free.
This leaves a precise design problem: to couple a low-dispersion physical
hierarchy with a scalable auxiliary inverse without allowing the latter to alter
the physical coarse correction.

We develop such a construction by combining an unshifted physical hierarchy
with a shifted auxiliary hierarchy.  One preconditioner application smooths on
the \(h\) grid and forms a physical correction with an independently
rediscretized IOFD operator \(A_{2h}\), rather than a Galerkin product.  A fixed
number of Krylov steps approximates \(A_{2h}^{-1}\), with each step
preconditioned by a complex-shifted \(2h\)--\(4h\) cycle.  The \(h\) and \(2h\)
equations therefore retain their propagative IOFD operators, while damping is
confined to the auxiliary machinery used to invert the coarse equation.  All
three levels are applied matrix-free, and the bottom correction is computed by
prescribed stencil and vector operations rather than a factorization.  The
resulting map has linear storage and right-preconditions flexible generalized
minimal residual (FGMRES) applied to the unshifted system \(A_h\vh=b_h\).

This construction yields a factorization-free IOFD hierarchy in which a target
resolution of six points per shortest wavelength places the \(2h\) correction at
only three.  Its viability therefore depends on both fine--coarse phase
compatibility and an effective auxiliary inverse.  LFA quantifies the former and
characterizes the contraction and finite-polynomial properties of the shifted
\(2h\)--\(4h\) hierarchy, thereby connecting the physical rediscretization and
the fixed auxiliary configuration within one three-grid analysis.

Numerical experiments separate wavefield validation from solver assessment.
The accuracy study validates the IOFD discretization against the outgoing Green
function.  The numerical wavefields reproduce both phase and relative amplitude,
with small full-domain relative errors.  The solver study then holds one
configuration fixed across the homogeneous, lens, wedge, barrier, and SEG/EAGE
Overthrust models, demonstrating robust convergence and scalable parallel
performance throughout.  In particular, a \revision{\(2560^3\)} problem spanning
approximately \revision{425} wavelengths per coordinate is solved in as little as
\revision{21.7 seconds} on 64 NVIDIA A100 GPUs.  CPU calculations
reach a \revision{\(13376^3\)}
problem spanning approximately \revision{2227} wavelengths per coordinate.  These experiments extend the scalability study to more than \revision{two trillion} unknowns.

The remainder of the paper is organized as follows.
Section~\ref{sec:iofd} formulates the PML-truncated problem and the IOFD
discretization.  Section~\ref{sec:three-grid} presents the nested
preconditioner, and Section~\ref{sec:three-grid-fourier} uses local Fourier
analysis to examine physical-grid phase compatibility and the shifted
auxiliary hierarchy.  Section~\ref{sec:validation} verifies
wavefield accuracy and establishes the solver configuration, while
Section~\ref{sec:scalability} evaluates convergence, memory use, and parallel
scalability on synthetic and geophysical models \revision{using GPUs and CPUs}.  Section~\ref{sec:conclusions}
draws the main conclusions and identifies directions for further work.

\section{Problem Formulation}
\label{sec:iofd}

This section specifies the PML-truncated Helmholtz problem and its IOFD
discretization.  We first formulate the radiating problem on a bounded domain
and then define the wavenumber-dependent 27-point operator and its amplitude
correction.

\subsection{PML-Truncated Problem}

Let \(\Omega_{\rm ph}\subset\mathbb{R}^3\) be a rectangular physical domain,
and let \(\bm{x}\) denote position.  We assume the time dependence
\(\exp(-\ii\omega t)\), where \(\ii^2=-1\) and \(\omega>0\) is the angular
frequency.  The acoustic pressure \(u\) satisfies
\begin{equation}
  -\Delta u(\bm{x})-k(\bm{x})^2u(\bm{x})=f(\bm{x}),
  \qquad
  k(\bm{x})=\frac{\omega}{c(\bm{x})},
  \label{eq:helmholtz-continuous}
\end{equation}
where \(c\) is the wave speed, \(k\) is the local wavenumber, and \(f\) is
compactly supported.  In a homogeneous exterior, let \(k_\infty\) denote the
constant exterior wavenumber and let \(r=|\bm{x}|\).  The outgoing solution is
selected by the Sommerfeld radiation condition
\begin{equation}
  \lim_{r\rightarrow\infty}r\bigl(\partial_r u-\ii k_\infty u\bigr)=0.
  \label{eq:sommerfeld}
\end{equation}
Here \(\partial_r\) denotes the radial derivative.

To truncate the exterior, \(\Omega_{\rm ph}\) is surrounded by a Cartesian
PML region.  For \(j\in\{1,2,3\}\), define
\begin{equation}
  \xi_j(x_j)=1+\ii\gamma_j(x_j),
  \qquad
  \gamma_j(x_j)=
  \gamma_{\max}\left[1-\cos\left(\frac{\pi d_j(x_j)}
  {2L_{\rm pml}}\right)\right],
  \label{eq:pml-profile}
\end{equation}
where \(d_j\in[0,L_{\rm pml}]\) is the distance measured into the PML in
coordinate direction \(j\) and vanishes in \(\Omega_{\rm ph}\), \(L_{\rm pml}\) is the
physical PML thickness, and \(\gamma_{\max}\) is the maximum stretch amplitude.
The resulting problem on the extended box \(\Omega\) is
\begin{equation}
  -\sum_{j=1}^3
  \xi_j^{-1}\partial_j
  \bigl(\xi_j^{-1}\partial_j u\bigr)
  -k(\bm{x})^2u=f
  \quad\text{in }\Omega,
  \qquad u=0\quad\text{on }\partial\Omega.
  \label{eq:helmholtz-pml}
\end{equation}
This coordinate-stretching formulation follows the PML construction used for
three-dimensional Helmholtz problems in \cite{Pinel2010}.  The wave speed is
extended into the PML by its nearest value in \(\Omega_{\rm ph}\); the
complex stretch then provides the attenuation in the exterior layer.

\subsection{IOFD Discretization}

Consider a uniform nodal grid with spacing \(h\) and grid points
\(\bm{x}_{\bm{i}}=h\bm{i}\), where \(\bm{i}\in\mathbb Z^3\) is a grid
multi-index.  In terms of the local number of points per wavelength (ppw), the
nondimensional wavenumber is
\begin{equation}
  \eta_{\bm{i}}=\frac{k(\bm{x}_{\bm{i}})h}{2\pi}
  =\frac{1}{\ppw(\bm{x}_{\bm{i}})}.
  \label{eq:local-inverse-ppw}
\end{equation}
The IOFD coefficient tables cover the interval \(0\leq\eta\leq0.4\).
Accordingly, for every level spacing
\(\ell\in\{h,2h,4h\}\), we define the interpolation coordinate
\begin{equation}
  \eta_{\ell,\bm i}^{\star}
  =\min\!\left\{\frac{k(\bm x_{\bm i})\ell}{2\pi},0.4\right\}.
  \label{eq:iofd-clamped-coordinate}
\end{equation}
All tabulated IOFD functions are evaluated at
\(\eta_{\ell,\bm i}^{\star}\), whereas the physical nondimensional
wavenumber \(\kappa_{\ell,\bm i}=k(\bm x_{\bm i})\ell\) is retained in the
mass factors.  Thus endpoint continuation of the coefficient functions is
part of every independently rediscretized level operator.
The IOFD scheme \cite{Stolk2016IOFD} separates the approximation into a
discrete Helmholtz operator \(P_h\), optimized for phase propagation, and a
compact operator \(Q_h\), chosen to correct the leading amplitude error.  The
five functions \(\alpha_j(s)\), \(j=1,\ldots,5\), \(s\in[0,0.4]\), defining
\(P_h\) are cubic Hermite interpolants of optimized tabulated values.
Consequently, the
27-point stencil adapts to the local value of \(kh\), rather than using one set
of coefficients throughout the domain.

Let
\(\mathcal{N}_m=\{\bm{\delta}\in\{-1,0,1\}^3:
\|\bm{\delta}\|_1=m\}\), \(m=0,\ldots,3\).  Away from the PML, the discrete
Helmholtz operator \(P_h\) has the shell form
\begin{equation}
  (P_hv)_{\bm{i}}=\frac{1}{h^2}
  \sum_{m=0}^3 C_m(\eta_{h,\bm{i}}^\star,\kappa_{\bm i})
  \sum_{\bm{\delta}\in\mathcal{N}_m}v_{\bm{i}+\bm{\delta}},
  \label{eq:iofd-shell-form}
\end{equation}
where, with \(a_j=\alpha_j(\eta_{h,\bm{i}}^\star)\) and
\(\kappa_{\bm{i}}=k(\bm{x}_{\bm{i}})h\),
\begin{align}
  C_0(\eta_{h,\bm i}^\star,\kappa_{\bm i}) &= 6a_4-\kappa_{\bm{i}}^2a_1, \nonumber\\
  C_1(\eta_{h,\bm i}^\star,\kappa_{\bm i}) &= -a_4+a_5-\frac{\kappa_{\bm{i}}^2}{6}a_2, \nonumber\\
  C_2(\eta_{h,\bm i}^\star,\kappa_{\bm i}) &= -\frac{a_5}{2}+\frac{1-a_4-a_5}{2}
         -\frac{\kappa_{\bm{i}}^2}{12}a_3, \nonumber\\
  C_3(\eta_{h,\bm i}^\star,\kappa_{\bm i}) &= -\frac{3}{4}(1-a_4-a_5)
         -\frac{\kappa_{\bm{i}}^2}{8}(1-a_1-a_2-a_3).
  \label{eq:iofd-shell-coefficients}
\end{align}
The four shells contain, respectively, one center point, six face neighbors,
twelve edge neighbors, and eight corner neighbors.  For variable wave speed,
the interpolation coordinate and physical wavenumber in a row are evaluated
at the central point of that row.

The tensor-product construction underlying
\eqref{eq:iofd-shell-form} also gives a direct extension to the PML\@.  In each
coordinate direction, the centered second difference is replaced by
\begin{equation}
  (D^{\rm pml}_{j,h}v)_{\bm{i}}
  =\frac{\xi_{j,\bm{i}}^{-1}}{h^2}
  \left[
  \xi_{j,\bm{i}+\frac12\bm{e}_j}^{-1}
  (v_{\bm{i}+\bm{e}_j}-v_{\bm{i}})
  -\xi_{j,\bm{i}-\frac12\bm{e}_j}^{-1}
  (v_{\bm{i}}-v_{\bm{i}-\bm{e}_j})
  \right].
  \label{eq:pml-second-difference}
\end{equation}
Here \(\bm e_j\) is the \(j\)th Cartesian unit vector.  The half-grid
stretching factors are obtained by averaging \(\gamma_j\) at
adjacent nodes before inversion.  Substituting
\eqref{eq:pml-second-difference} for each directional second difference in
the IOFD tensor product preserves the 27-point support.  The resulting stencil
has transverse weights \(a_4\), \(a_5/4\), and
\((1-a_4-a_5)/4\) at the transverse center, edge, and corner positions,
respectively, while the mass weights are \(a_1\), \(a_2/6\), \(a_3/12\), and
\((1-a_1-a_2-a_3)/8\) on the four three-dimensional shells.  Thus the compact
structure is retained, although directional stretching removes the shell
symmetry inside the PML\@.

IOFD also introduces a compact amplitude-correction operator \(Q_h\).  Let
\(\beta_j(s)\), \(j=1,2,3\), \(s\in[0,0.4]\), denote the corresponding
Hermite interpolants and define
\begin{equation}
\begin{aligned}
  \mathcal Q_0(s)&=\beta_1(s), &
  \mathcal Q_1(s)&=\frac{\beta_2(s)}{6},\\
  \mathcal Q_2(s)&=\frac{\beta_3(s)}{12}, &
  \mathcal Q_3(s)&=
  \frac{1-\beta_1(s)-\beta_2(s)-\beta_3(s)}{8}.
\end{aligned}
  \label{eq:q-shell-coefficients}
\end{equation}
Then
\begin{equation}
  (Q_hw)_{\bm{i}}=
  \sum_{m=0}^3 \mathcal Q_m(\eta_{h,\bm{i}}^\star)
  \sum_{\bm{\delta}\in\mathcal{N}_m}w_{\bm{i}+\bm{\delta}}.
  \label{eq:q-operator}
\end{equation}
Let \(f_h\) denote the nodal discretization of \(f\), and let \(v_h\) be the
auxiliary grid field used by IOFD\@.  The complete discretization and the
reconstructed pressure \(u_h\) are therefore
\begin{equation}
  P_h\vh=Q_hf_h,\qquad \uh=Q_h\vh.
  \label{eq:iofd-complete-system}
\end{equation}
The first application of \(Q_h\) corrects the discrete representation of the
source, whereas the second reconstructs the pressure field with the IOFD
amplitude correction.  In the remainder of the paper we set
\(A_h=P_h\) and \(b_h=Q_hf_h\).  Krylov iteration is applied to
\(A_h\vh=b_h\), and the reported relative residual is measured as
\begin{equation}
  \frac{\|A_h\vh-b_h\|_2}{\|b_h\|_2}.
  \label{eq:relative-residual}
\end{equation}
The operator \(Q_h\) is not part of an individual preconditioner application.

\section{Three-Grid Preconditioner}
\label{sec:three-grid}

This section develops the three-grid preconditioner for \(A_h\vh=b_h\).  We
first define the independently rediscretized IOFD operators and geometric grid
transfers, then construct the shifted \(2h\)-to-\(4h\) auxiliary cycle and the
complete flexible preconditioner.  We retain unbold symbols for the global grid
fields and use boldface for temporary vectors within a cycle.

\subsection{Physical and Auxiliary Operators}

Let \(\mathcal X_h\), \(\mathcal X_{2h}\), and \(\mathcal X_{4h}\) denote the
spaces of complex-valued grid functions on three nested nodal grids.  Between
consecutive levels, the mesh spacing is doubled and the number of PML layers is
halved, so the physical PML thickness is retained.  On each level
\(\ell\in\{h,2h,4h\}\), let \(L_\ell\) denote the compact discretization of
the negative PML-stretched Laplacian and let \(M_\ell\) denote the IOFD weighted discretization of
\(k(\bm{x})^2u\).  Both operators are independently rediscretized using the
level spacing \(\ell\), the wave speed sampled at the corresponding grid
points, and the interpolation coordinate
\(\eta_{\ell,\bm i}^{\star}\) in
\eqref{eq:iofd-clamped-coordinate}.  On the two physical levels, define
\begin{equation}
  A_\ell=L_\ell-M_\ell,\qquad \ell\in\{h,2h\},
  \label{eq:physical-iofd-levels}
\end{equation}
using the construction in Section~\ref{sec:iofd}.  In particular,
\begin{equation}
  A_{2h}\ne\mathcal R_h A_h\mathcal P_h.
  \label{eq:not-galerkin}
\end{equation}
Instead, the coefficients of \(A_{2h}\) are evaluated at the local coordinate
\(\eta_{2h,\bm i}^{\star}\), while the mass factors retain
\(k(\bm{x}_{\bm i})\,2h\).  This distinction is essential because the IOFD
coarse operator is designed to reproduce wave propagation at the coarse
resolution rather than to inherit the fine-grid stencil algebraically
\cite{StolkAhmedBhowmik2014,Stolk2016IOFD}.

The prolongations
\(\mathcal P_h:\mathcal X_{2h}\rightarrow\mathcal X_h\) and
\(\mathcal P_{2h}:\mathcal X_{4h}\rightarrow\mathcal X_{2h}\) are trilinear.
In one coordinate direction their action is
\begin{equation}
  (\mathcal P_\ell w)_{2j}=w_j,\qquad
  (\mathcal P_\ell w)_{2j+1}=\frac12(w_j+w_{j+1}).
  \label{eq:linear-prolongation}
\end{equation}
The restriction uses full weighting.  In one coordinate direction,
\begin{equation}
  (\mathcal R_\ell v)_j
  =\frac14v_{2j-1}+\frac12v_{2j}+\frac14v_{2j+1}.
  \label{eq:full-weighting}
\end{equation}
Consequently, in three dimensions
\(\mathcal R_\ell=2^{-3}\mathcal P_\ell^*\), where \(^*\) denotes the
conjugate transpose.

The shifted hierarchy is introduced only on the two coarser grids.  For
\(\ell\in\{2h,4h\}\), define
\begin{equation}
  B_\ell=L_\ell-(1+\ii\sigma)M_\ell,\qquad \sigma>0.
  \label{eq:shifted-iofd-levels}
\end{equation}
The shift multiplies the physical mass factor
\(\kappa_{\ell,\bm i}^2\); it does not change the real interpolation
coordinate \(\eta_{\ell,\bm i}^{\star}\) used to evaluate the IOFD tables.
The complex shift damps components that are poorly represented on the coarse
grids and thereby makes standard relaxation effective.  It does not alter
either the target equation \(A_h\vh=b_h\) or the unshifted intermediate
equation involving \(A_{2h}\).

\subsection{Shifted Auxiliary Cycle}

The auxiliary cycle acts only on the two shifted coarse levels.  We first define
the relaxation used on these levels.  For an operator \(H\) with diagonal
\(D_H\), a relaxation weight \(\omega_J\), and an integer \(n_J\geq1\), let
\(\mathcal J_{n_J}(H,\bm g;\omega_J)\) denote \(n_J\) damped Jacobi sweeps
initialized from zero,
\begin{equation}
\begin{aligned}
  \bm y^{(0)}&=\bm 0,\\
  \bm y^{(j+1)}
  &=\bm y^{(j)}
  +\omega_JD_H^{-1}\bigl(\bm g-H\bm y^{(j)}\bigr),
  &&j=0,\ldots,n_J-1,\\
  \mathcal J_{n_J}(H,\bm g;\omega_J)&=\bm y^{(n_J)}.
\end{aligned}
  \label{eq:jacobi-map}
\end{equation}
Let \(n_{J,\ell}\) denote the Jacobi sweep count used on level \(\ell\).  At the
bottom level, \(\mathcal C_{4h}(\bm g)\) denotes \(m_{4h}\) prescribed
communication-avoiding generalized minimal residual (CA-GMRES) steps for
\(B_{4h}\bm e=\bm g\), initialized from zero and right-preconditioned by
\(\mathcal J_{n_{J,4h}}(B_{4h},\cdot;\omega_{4h})\).
Here and below, the subscripted quantities \(\omega_\ell\) are relaxation
weights and are unrelated to the angular frequency \(\omega\) in
\eqref{eq:helmholtz-continuous}.

One shifted \(2h\)-to-\(4h\) V-cycle defines the map
\(\mathcal V_{2h}^{\sigma}:\mathcal X_{2h}\rightarrow\mathcal X_{2h}\).
For a residual \(\bm q_{2h}\), it is given by
\begin{equation}
\begin{aligned}
  \bm y_{2h}
  &=\mathcal J_{n_{J,2h}}(B_{2h},\bm q_{2h};\omega_{2h}),\\
  \bm q_{4h}
  &=\mathcal R_{2h}
    \bigl(\bm q_{2h}-B_{2h}\bm y_{2h}\bigr),\\
  \bm e_{4h}
  &=\mathcal C_{4h}(\bm q_{4h}),\\
  \widetilde{\bm y}_{2h}
  &=\bm y_{2h}+\mathcal P_{2h}\bm e_{4h},\\
  \mathcal V_{2h}^{\sigma}(\bm q_{2h})
  &=\widetilde{\bm y}_{2h}
    +\mathcal J_{n_{J,2h}}\!\left(
       B_{2h},\bm q_{2h}-B_{2h}\widetilde{\bm y}_{2h};
       \omega_{2h}\right).
\end{aligned}
  \label{eq:shifted-auxiliary-cycle}
\end{equation}
Thus the \(4h\) problem is treated by prescribed Krylov work rather than a
direct factorization.  The unshifted intermediate correction is then defined
by
\begin{equation}
  \bm e_{2h}
  =\mathcal K_{2h}(\bm r_{2h}),
  \label{eq:unshifted-middle-solve}
\end{equation}
where \(\mathcal K_{2h}\) consists of \(c_{2h}\) restarted FGMRES cycles with
restart length \(m_{2h}\) for \(A_{2h}\bm e_{2h}=\bm r_{2h}\), initialized from zero and
right-preconditioned by \(\mathcal V_{2h}^{\sigma}\).  Hence the operator
driving this Krylov process is the propagative, unshifted \(A_{2h}\), while
the shift appears only in the approximate inverse used within it.

\subsection{Three-Grid Cycle}

We now assemble the physical correction and the shifted auxiliary solve.  Let
\(\mathcal S_h(\bm g;\bm x_0)\) denote one right-preconditioned GMRES cycle of
length \(m_h\) for \(A_h\bm x=\bm g\), initialized by \(\bm x_0\).  Its right
preconditioner is \(\mathcal J_{n_{J,h}}(A_h,\cdot;\omega_h)\).
Algorithm~\ref{alg:three-grid} then gives one application of
\(\mathcal M_h^{-1}\).

\begin{algorithm}[htbp]
\caption{Application of the three-grid preconditioner
\(\bm z_h=\mathcal M_h^{-1}\bm r_h\).}
\label{alg:three-grid}
\begin{enumerate}
  \renewcommand{\labelenumi}{\arabic{enumi}:}
  \setlength{\itemsep}{0.12em}
  \setlength{\parsep}{0pt}
  \setlength{\parskip}{0pt}
  \item Pre-smooth on the fine grid:
  \(\bm z_h^{\rm pre}=\mathcal S_h(\bm r_h;\bm 0)\).
  \item Restrict the fine-grid residual:
  \(\bm r_{2h}=\mathcal R_h
  (\bm r_h-A_h\bm z_h^{\rm pre})\).
  \item Compute the unshifted intermediate-grid correction:
  \(\bm e_{2h}=\mathcal K_{2h}(\bm r_{2h})\).
  \item Prolong and correct:
  \(\widetilde{\bm z}_h
  =\bm z_h^{\rm pre}+\mathcal P_h\bm e_{2h}\).
  \item Post-smooth on the fine grid:
  \(\bm z_h=\mathcal S_h(\bm r_h;\widetilde{\bm z}_h)\).
\end{enumerate}
\end{algorithm}

The inner Krylov maps depend on the right-hand sides from which their Krylov
spaces are generated, so \(\mathcal M_h^{-1}\) is generally nonlinear.  The
outer iteration is FGMRES with restart length \(m_{\rm out}\)
\cite{Saad1993FGMRES}, using \(\mathcal M_h^{-1}\) as a right preconditioner for
\(A_h\vh=b_h\).
Only the outer iteration performs a convergence test, based on
\eqref{eq:relative-residual}; every inner method executes its prescribed number of
steps.  Each outer Arnoldi step invokes the preconditioner once, and the
reported number of preconditioner calls is therefore the number of outer
FGMRES steps.  The integer parameters \(n_{J,\ell}\), \(m_h\), \(c_{2h}\),
\(m_{2h}\), \(m_{4h}\), and \(m_{\rm out}\) are specified with the numerical
configuration in Section~\ref{sec:scalability}.

\section{Local Fourier Analysis}
\label{sec:three-grid-fourier}

The hierarchy requires a phase-compatible physical correction and an effective
auxiliary inverse.  We examine these properties through the fine--coarse
dispersion relation, the Fourier blocks of the shifted cycle, and the residual
polynomials of the bottom solver.  Together, these diagnostics provide a
source-independent basis for choosing the auxiliary parameters.

\subsection{Local Wavelength Range}

The fine-grid spacing is chosen so that the shortest wavelength, attained at
\(c_{\min}\), has \(\ppw_{\min}=6\) points.  If
\(\chi=c_{\max}/c_{\min}\) is the velocity contrast in the propagating medium,
then the local fine-grid sampling is
\begin{equation}
  \nu(\bm x):=\ppw(\bm x)
  =\frac{2\pi}{k(\bm x)h}
  =\ppw_{\min}\frac{c(\bm x)}{c_{\min}}.
  \label{eq:local-ppw-range}
\end{equation}
We consider the practically relevant contrast range \(\chi\leq4\), and hence
\(\nu(\bm x)\in[6,24]\).  This range covers, for example, the velocity ranges
represented by the SEG/EAGE three-dimensional
Salt and Overthrust benchmarks \cite{AminzadehBracKunz1997}.  The sampling
interval is therefore fixed without reference to a source or to a particular
realization of the velocity field.
The angular and local-wavelength scans used below are refined until the
reported digits are unchanged.

\subsection{Phase Compatibility of the Physical Grids}

For a locally homogeneous interior stencil, denote the Fourier symbol of
\(A_\ell\), \(\ell\in\{h,2h\}\), by
\(\widehat a_\ell(\bm\theta;\nu)\).  Let
\(\mathcal D_r\subset\mathbb S^2\) and
\(\mathcal N_r\subset[6,24]\) denote the nested directional and local-sampling
sets at refinement level \(r\).  For each
\((\bm d,\nu)\in\mathcal D_r\times\mathcal N_r\), let
\(q_\ell(\bm d;\nu)\) be the normalized numerical wavenumber on the physical
branch:
\begin{equation}
  \widehat a_\ell\bigl(q_\ell(\bm d;\nu)\,k\ell\bm d;\nu\bigr)=0.
  \label{eq:normalized-numerical-wavenumber}
\end{equation}
Here \(k=2\pi/(\nu h)\), and \(q_\ell=1\) corresponds to the exact physical
wavenumber.  At every sample the real root is bracketed in \([0.9,1.1]\) and
computed by 50 bisection steps; the same bracket remains valid throughout the
refined sets.  We write \(\ppw_\ell=2\pi/(k\ell)\) for the points per wavelength on
level \(\ell\).  On the finest sets, the largest sampled fine--coarse mismatch
occurs at the shortest wavelength, \(\nu=6\), where
\begin{equation}
\begin{aligned}
  \max_{\bm d\in\mathcal D_r}|q_h(\bm d;6)-1|
  &=2.96\times10^{-6},\\
  \max_{\bm d\in\mathcal D_r}|q_{2h}(\bm d;6)-1|
  &=3.27\times10^{-4},\\
  \max_{\bm d\in\mathcal D_r}|q_{2h}(\bm d;6)-q_h(\bm d;6)|
  &=3.24\times10^{-4},
\end{aligned}
  \qquad
  \ppw_h=\nu=6,\quad \ppw_{2h}=\nu/2=3.
  \label{eq:iofd-phase-matching}
\end{equation}
Thus, over the refined angular and local-wavelength sets, the first coarse grid reproduces
the normalized fine-grid numerical wavenumber to within
\(3.3\times10^{-4}\), even though it contains only three points per shortest
wavelength.  The sampled mismatch decreases monotonically over
\(\mathcal N_r\).  This quantitative phase evidence motivates rediscretizing
the physical coarse operator with IOFD: a generic Galerkin product is not
constructed to preserve the same dispersion relation.

\subsection{Fourier Blocks of the Auxiliary Hierarchy}

The auxiliary hierarchy is analyzed locally on an infinite grid after freezing
the interior coefficients.  This standard LFA setting gives the
translation-invariant symbols relevant to smoothing and coarse correction.  The
stationary discretization, transfer, and Jacobi components are analyzed directly
within the multilevel harmonic framework
\cite{WienandsOosterlee2001,CoolsVanroose2013,CalandraEtAl2013}.

Let \(\kappa_\ell=k\ell=2\pi(\ell/h)/\nu\).  Under this local model,
\eqref{eq:iofd-clamped-coordinate} reduces to
\[
  \eta_\ell^\star=\min\!\left\{\frac{\kappa_\ell}{2\pi},0.4\right\}.
\]
Thus the symbol below is
the interior symbol of the same rediscretized level operator used in
the algorithm: the IOFD tables are evaluated at \(\eta_\ell^\star\), while
the physical value of \(\kappa_\ell\) is retained in the mass term.  Let
\(C_{m,\ell}^{\sigma}(\nu)\), \(m=0,\ldots,3\), denote the four
shell coefficients in \eqref{eq:iofd-shell-coefficients}, evaluated at
\(\eta_\ell^\star\), with \(\kappa_\ell^2\) replaced by
\((1+\ii\sigma)\kappa_\ell^2\).  The symbol of the shifted operator is then
\begin{equation}
\begin{aligned}
  \widehat b_\ell^{\sigma}(\bm\theta;\nu)
  =\ell^{-2}\bigl[&C_{0,\ell}^{\sigma}
  +2C_{1,\ell}^{\sigma}(\zeta_1+\zeta_2+\zeta_3)\\
  &+4C_{2,\ell}^{\sigma}(\zeta_1\zeta_2+\zeta_1\zeta_3+\zeta_2\zeta_3)
  +8C_{3,\ell}^{\sigma}\zeta_1\zeta_2\zeta_3\bigr],
  \qquad \zeta_j=\cos\theta_j.
\end{aligned}
  \label{eq:shifted-iofd-symbol}
\end{equation}
Thus \(\widehat a_\ell=\widehat b_\ell^0\) is the symbol of the unshifted
IOFD operator.  Unlike the standard second-order Helmholtz symbol, both the
stiffness distribution and the mass distribution in
\eqref{eq:shifted-iofd-symbol} depend on \(k\ell\).
We suppress the dependence of these symbols on \(\nu\) through
\(\kappa_\ell\) whenever it is unambiguous.

The dimensionless mass symbol is
\[
\widehat m(\bm\theta;\eta)
=a_1+\frac{a_2}{3}\sum_{j=1}^3\zeta_j
+\frac{a_3}{3}\sum_{1\leq i<j\leq3}\zeta_i\zeta_j
+(1-a_1-a_2-a_3)\zeta_1\zeta_2\zeta_3,
\qquad a_j=\alpha_j(\eta).
\]
The inverse symbols used below are well defined for
\(\sigma\in[0.75,1]\).  Evaluating the cosine vertices and cubic-Hermite
stationary points of the dimensionless mass symbol
\(\widehat m(\bm\theta;\eta)\) gives
\[
  \widehat m(\bm\theta;\eta)\geq0.489,\qquad
  \alpha_1(\eta)\geq0.588,\qquad
  0\leq\eta\leq0.4.
\]
Thus
\(\operatorname{Im}(\ell^2\widehat b_\ell^\sigma)
=-\sigma\kappa_\ell^2\widehat m\ne0\) and
\(\operatorname{Im}(\ell^2d_\ell^\sigma)
=-\sigma\kappa_\ell^2a_1\ne0\) on the shifted levels.  Direct evaluation also
gives \(h^2d_h^0\geq3.519\) for \(\nu\in[6,24]\), proving that every symbol
or Jacobi diagonal inverted below is nonzero.

Consider first the shifted \(2h\)-to-\(4h\) cycle.  For
\(\bm\theta\in\Theta_{\rm low}=(-\pi/2,\pi/2]^3\), define the eight
\(4h\)-harmonics on the \(2h\) grid by
\begin{equation}
  \bm\theta^{\bm p}=\bm\theta+\pi\bm p,
  \qquad \bm p\in\{0,1\}^3.
  \label{eq:two-grid-harmonics}
\end{equation}
The angles are understood modulo \(2\pi\).
The associated harmonic space is invariant under the translation-invariant
IOFD stencil and weighted Jacobi iteration.  In this basis, set
\begin{equation}
  \widehat{\bm B}_{2h}^{\sigma}(\bm\theta;\nu)
  =\operatorname{diag}_{\bm p}
   \widehat b_{2h}^{\sigma}(\bm\theta^{\bm p};\nu),
  \qquad
  \widehat{\bm A}_{2h}(\bm\theta;\nu)
  =\widehat{\bm B}_{2h}^{0}(\bm\theta;\nu).
  \label{eq:lfa-block-operators}
\end{equation}
We represent Fourier coefficients in the unnormalized harmonic bases
\[
  \phi_{\bm p}(\bm n)
  =\exp\!\bigl(\ii(\bm\theta+\pi\bm p)\cdot\bm n\bigr),
  \qquad
  \psi(\bm j)=\exp(2\ii\bm\theta\cdot\bm j).
\]
In these bases, trilinear interpolation and full weighting are represented by
the column \(\widehat{\bm p}(\bm\theta)\) and its adjoint, whose entries are
\begin{equation}
  \widehat p_{\bm p}(\bm\theta)
  =\prod_{j=1}^3\frac{1+\cos(\theta_j+\pi p_j)}{2}.
  \label{eq:lfa-transfer-symbol}
\end{equation}
Direct application of \eqref{eq:linear-prolongation} and
\eqref{eq:full-weighting}, followed by tensorization, gives
\[
  \mathcal P_{2h}\psi
  =\sum_{\bm p}\widehat p_{\bm p}\phi_{\bm p},
  \qquad
  \mathcal R_{2h}\phi_{\bm p}
  =\widehat p_{\bm p}\psi.
\]
Hence the coefficient matrices are \(\widehat{\bm p}\) and
\(\widehat{\bm p}^{*}\).  This is consistent with
\(\mathcal R_{2h}=2^{-3}\mathcal P_{2h}^{*}\): over a periodic cell,
\(\|\phi_{\bm p}\|_2^2=2^3\|\psi\|_2^2\), and the basis-norm ratio supplies
the factor \(2^{-3}\) \cite{WienandsOosterlee2001}.
Here and below, \(\bm I\) denotes an identity matrix of the appropriate
dimension.

Let \(d_{\ell}^{\sigma}(\nu)=C_{0,\ell}^{\sigma}(\nu)/\ell^2\).  On the
\(2h\) harmonic space, one Jacobi error sweep has the diagonal block
\begin{equation}
  \widehat{\bm S}_{2h}^{\sigma}
  (\bm\theta;\nu,\omega_{2h})
  =\bm I-\omega_{2h}\bigl(d_{2h}^{\sigma}(\nu)\bigr)^{-1}
  \widehat{\bm B}_{2h}^{\sigma}(\bm\theta;\nu).
  \label{eq:lfa-jacobi-symbol}
\end{equation}
For \(\ell\in\{h,2h,4h\}\), with \(\sigma=0\) on the fine grid, define
\(\Theta_{\rm high}=(-\pi,\pi]^3\setminus\Theta_{\rm low}\).  After
\(n_{J,\ell}\) sweeps, the high-frequency smoothing factor at local fine-grid
sampling \(\nu\) is
\begin{equation}
\begin{aligned}
  \mu_{\ell}(\nu;\sigma,\omega_{\ell},n_{J,\ell})
  &=
  \sup_{\bm\vartheta\in\Theta_{\rm high}}
  \left|1-\omega_{\ell}
  \frac{\widehat b_{\ell}^{\sigma}(\bm\vartheta;\nu)}
       {d_{\ell}^{\sigma}(\nu)}\right|^{n_{J,\ell}},
  \\
  \mu_{\ell}^{\rm wc}(\sigma,\omega_{\ell},n_{J,\ell})
  &=
  \sup_{\nu\in[6,24]}
  \mu_{\ell}(\nu;\sigma,\omega_{\ell},n_{J,\ell}).
\end{aligned}
  \label{eq:lfa-smoothing-factor}
\end{equation}
Here the superscript \({\rm wc}\) denotes the worst case over
\(\nu\in[6,24]\).
For fixed \(\nu\), the Jacobi factor is affine in each
\(\zeta_j=\cos\vartheta_j\), and its modulus is separately convex.  The cosine
image of \(\Theta_{\rm high}\) is a union of seven boxes with vertices among
\(\{-1,0,1\}^3\setminus\{(1,1,1)\}\).  Therefore these 26 candidates give the
angular supremum in \eqref{eq:lfa-smoothing-factor}; only the remaining
one-dimensional supremum over \(\nu\) requires refinement.
To isolate the effects of the shift and relaxation, we first take an exact
\(4h\) solve.  The resulting reference cycle has the error block
\begin{equation}
\begin{aligned}
  \widehat{\bm E}_{2h}^{\sigma}
  (\bm\theta;\nu,\omega_{2h},n_{J,2h})
  ={}&(\widehat{\bm S}_{2h}^{\sigma})^{n_{J,2h}}
  \left[\bm I-
  \widehat{\bm p}
  \frac{1}{\widehat b_{4h}^{\sigma}(2\bm\theta;\nu)}
  \widehat{\bm p}^{*}
  \widehat{\bm B}_{2h}^{\sigma}\right]
  (\widehat{\bm S}_{2h}^{\sigma})^{n_{J,2h}}.
\end{aligned}
  \label{eq:lfa-shifted-vcycle}
\end{equation}
The repeated arguments \((\bm\theta,\nu,\omega_{2h})\) have been suppressed on
the right-hand side.
Consequently, one shifted cycle initialized from zero has the approximate
inverse symbol
\begin{equation}
  \widehat{\bm V}_{2h}^{\sigma}
  (\bm\theta;\nu,\omega_{2h},n_{J,2h})
  =\left(\bm I-\widehat{\bm E}_{2h}^{\sigma}\right)
   \left(\widehat{\bm B}_{2h}^{\sigma}\right)^{-1}.
  \label{eq:lfa-shifted-inverse}
\end{equation}
The associated preconditioned coarse blocks are the eight-dimensional matrices
\begin{equation}
  \widehat{\bm H}_{2h}^{\sigma}
  (\bm\theta;\nu,\omega_{2h},n_{J,2h})
  =\widehat{\bm A}_{2h}(\bm\theta;\nu)
   \widehat{\bm V}_{2h}^{\sigma}
   (\bm\theta;\nu,\omega_{2h},n_{J,2h}).
  \label{eq:lfa-preconditioned-coarse}
\end{equation}
In the following identities, all arguments are again suppressed.  Let
\(\widehat{\bm M}_{2h}\) be the harmonic block of the mass operator.  Since
\(\widehat{\bm B}_{2h}^{\sigma}
=\widehat{\bm A}_{2h}-\ii\sigma\widehat{\bm M}_{2h}\),
\begin{equation}
\begin{aligned}
  \widehat{\bm H}_{2h}^{\sigma}
  &=\widehat{\bm A}_{2h}
    (\bm I-\widehat{\bm E}_{2h}^{\sigma})
    (\widehat{\bm B}_{2h}^{\sigma})^{-1},\\
  \widehat{\bm H}_{2h}^{\sigma}
  -\widehat{\bm A}_{2h}(\widehat{\bm B}_{2h}^{\sigma})^{-1}
  &=-\widehat{\bm A}_{2h}\widehat{\bm E}_{2h}^{\sigma}
    (\widehat{\bm B}_{2h}^{\sigma})^{-1},\\
  \widehat{\bm A}_{2h}(\widehat{\bm B}_{2h}^{\sigma})^{-1}
  &=\bm I+\ii\sigma\widehat{\bm M}_{2h}
    (\widehat{\bm B}_{2h}^{\sigma})^{-1}.
\end{aligned}
  \label{eq:lfa-shift-balance}
\end{equation}
Thus the quality of the auxiliary inverse is governed jointly by the
multigrid approximation error and the shift-induced perturbation.  Their
separation in \eqref{eq:lfa-shift-balance} identifies the competing effects
to be balanced by the auxiliary parameters.
This construction combines the Fourier representation of a multigrid
approximate inverse \cite{WienandsOosterleeWashio2000} with the spectral
analysis of shifted-Laplacian preconditioning
\cite{VanGijzenErlanggaVuik2007}; here the IOFD symbol replaces the standard
finite-difference symbol.

The same representation also explains the bottom solver.  Define the scalar
Jacobi error factor
\begin{equation}
  \widehat s_{4h}^{\sigma}
  (\bm\theta;\nu,\omega_{4h})
  =1-\omega_{4h}
   \frac{\widehat b_{4h}^{\sigma}(\bm\theta;\nu)}
        {d_{4h}^{\sigma}(\nu)}.
  \label{eq:lfa-bottom-jacobi-factor}
\end{equation}
Then \(n_{J,4h}\) Jacobi sweeps initialized from zero define the polynomial
approximate inverse
\begin{equation}
  \widehat\Xi_{4h}^{\sigma}
  (\bm\theta;\nu,\omega_{4h},n_{J,4h})
  =\frac{\omega_{4h}}{d_{4h}^{\sigma}(\nu)}
  \sum_{j=0}^{n_{J,4h}-1}
  \bigl(\widehat s_{4h}^{\sigma}
  (\bm\theta;\nu,\omega_{4h})\bigr)^j.
  \label{eq:lfa-bottom-inverse}
\end{equation}
Since
\(1-\widehat s_{4h}^{\sigma}
=\omega_{4h}\widehat b_{4h}^{\sigma}/d_{4h}^{\sigma}\), telescoping gives
the right-preconditioned symbol
\begin{equation}
  \widehat b_{4h}^{\sigma}\widehat\Xi_{4h}^{\sigma}
  =1-(\widehat s_{4h}^{\sigma})^{n_{J,4h}}.
  \label{eq:lfa-bottom-preconditioned-symbol}
\end{equation}
Whenever \(\widehat b_{4h}^{\sigma}\ne0\),
\eqref{eq:lfa-bottom-inverse} is equivalently
\([1-(\widehat s_{4h}^{\sigma})^{n_{J,4h}}]
(\widehat b_{4h}^{\sigma})^{-1}\).

\subsection{Parameter Diagnostics}

For the reference cycle in \eqref{eq:lfa-shifted-vcycle}, define the
contraction factor
\begin{equation}
  q_V(\sigma,\omega_{2h};n_{J,2h})
  =
  \sup_{\substack{\nu\in[6,24]\\
                  \bm\theta\in\Theta_{\rm low}}}
  \left\|
    \widehat{\bm E}_{2h}^{\sigma}
    (\bm\theta;\nu,\omega_{2h},n_{J,2h})
  \right\|_2.
  \label{eq:worst-vcycle-factor}
\end{equation}
The Fourier transform block-diagonalizes this cycle; hence \(q_V<1\) is
sufficient for one such cycle to contract in the discrete
\(\ell_2\) norm for every local sampling in \([6,24]\).  The matrix two-norm is
essential here: the harmonic blocks need not be normal, so their spectral
radii do not bound one-cycle error amplification.
The factor \(q_V\) isolates the shift and \(2h\) relaxation.  We pair it below
with an \(m_{4h}\)-dependent diagnostic for the finite bottom work.

The bottom correction \(\mathcal C_{4h}\) uses \(m_{4h}\) CA-GMRES steps,
whose residual polynomials quantify the prescribed bottom work.  For its
\(n_{J,4h}\)-sweep Jacobi right preconditioner, let
\begin{equation}
\begin{aligned}
  \widehat\lambda_{4h}^{\sigma}
  (\bm\theta;\nu,\omega_{4h},n_{J,4h})
  &=1-\bigl(\widehat s_{4h}^{\sigma}
      (\bm\theta;\nu,\omega_{4h})\bigr)^{n_{J,4h}},
  \\
  \rho_{4h}(\sigma,\omega_{4h};n_{J,4h})
  &=\sup_{\nu,\bm\theta}|\widehat\lambda_{4h}^{\sigma}|,\\
  \delta_{4h}(\sigma,\omega_{4h};n_{J,4h})
  &=\inf_{\nu,\bm\theta}|\widehat\lambda_{4h}^{\sigma}|.
\end{aligned}
  \label{eq:bottom-radius-indicator}
\end{equation}
Here both extrema are taken over \(\nu\in[6,24]\) and
\(\bm\theta\in(-\pi,\pi]^3\).  For each fixed \(\nu\), the local Fourier
multiplier is normal and its spectrum is the range of
\(\widehat\lambda_{4h}^{\sigma}\).  The extrema in
\eqref{eq:bottom-radius-indicator} therefore provide outer-radius and
origin-separation diagnostics over the local sampling interval.  These
characterize nonsingularity and spectral scaling.  Let \(\Pi_m^{\mathbb{C}}\)
denote the space of complex polynomials of degree at most \(m\).  Finite bottom
work is quantified by the complementary degree-dependent ideal residual factor
\begin{equation}
  \gamma_{4h}^{(m)}
  (\sigma,\omega_{4h};n_{J,4h})
  =
  \inf_{\substack{p\in\Pi_m^{\mathbb{C}}\\p(0)=1}}
  \sup_{\substack{\nu\in[6,24]\\
                  \bm\theta\in(-\pi,\pi]^3}}
  \left|p\!\left(\widehat\lambda_{4h}^{\sigma}
  (\bm\theta;\nu,\omega_{4h},n_{J,4h})\right)\right|.
  \label{eq:bottom-ideal-gmres-factor}
\end{equation}
The factor \(\gamma_{4h}^{(m)}\) bounds the residual reduction of the
bottom solve in the local Fourier model.  To see this, let
\(H_{4h}^{\sigma}(\nu)=B_{4h}^{\sigma}(\nu)\Xi_{4h}^{\sigma}(\nu)\)
be a periodic realization of the right-preconditioned operator, where
\(\Xi_{4h}^{\sigma}(\nu)\) is the Jacobi approximate inverse with symbol
\eqref{eq:lfa-bottom-inverse}.  The product has symbol
\(\widehat\lambda_{4h}^{\sigma}\).  The unitary discrete Fourier transform
diagonalizes \(H_{4h}^{\sigma}(\nu)\), so the norm of any polynomial in this
operator equals the largest modulus of that polynomial on its Fourier
eigenvalues.  GMRES residual minimization therefore gives
\cite{SaadSchultz1986}
\begin{equation}
  \frac{\|\bm r_m\|_2}{\|\bm r_0\|_2}
  \leq
  \inf_{\substack{p\in\Pi_m^{\mathbb C}\\p(0)=1}}
  \bigl\|p\bigl(H_{4h}^{\sigma}(\nu)\bigr)\bigr\|_2
  \leq \gamma_{4h}^{(m)},
  \qquad \nu\in[6,24],
  \label{eq:bottom-gmres-residual-bound}
\end{equation}
where \(\bm r_j\) is the bottom-grid residual after \(j\) steps and
\(\bm r_0\ne\bm0\).  The final inequality follows because
\eqref{eq:bottom-ideal-gmres-factor} uses one polynomial over the entire
local-sampling interval.  Taking \(m=m_{4h}\) thus connects the bottom
iteration count to a residual bound uniform in the right-hand side.

\begin{figure}[t]
  \centering
  \includegraphics[width=0.82\textwidth]{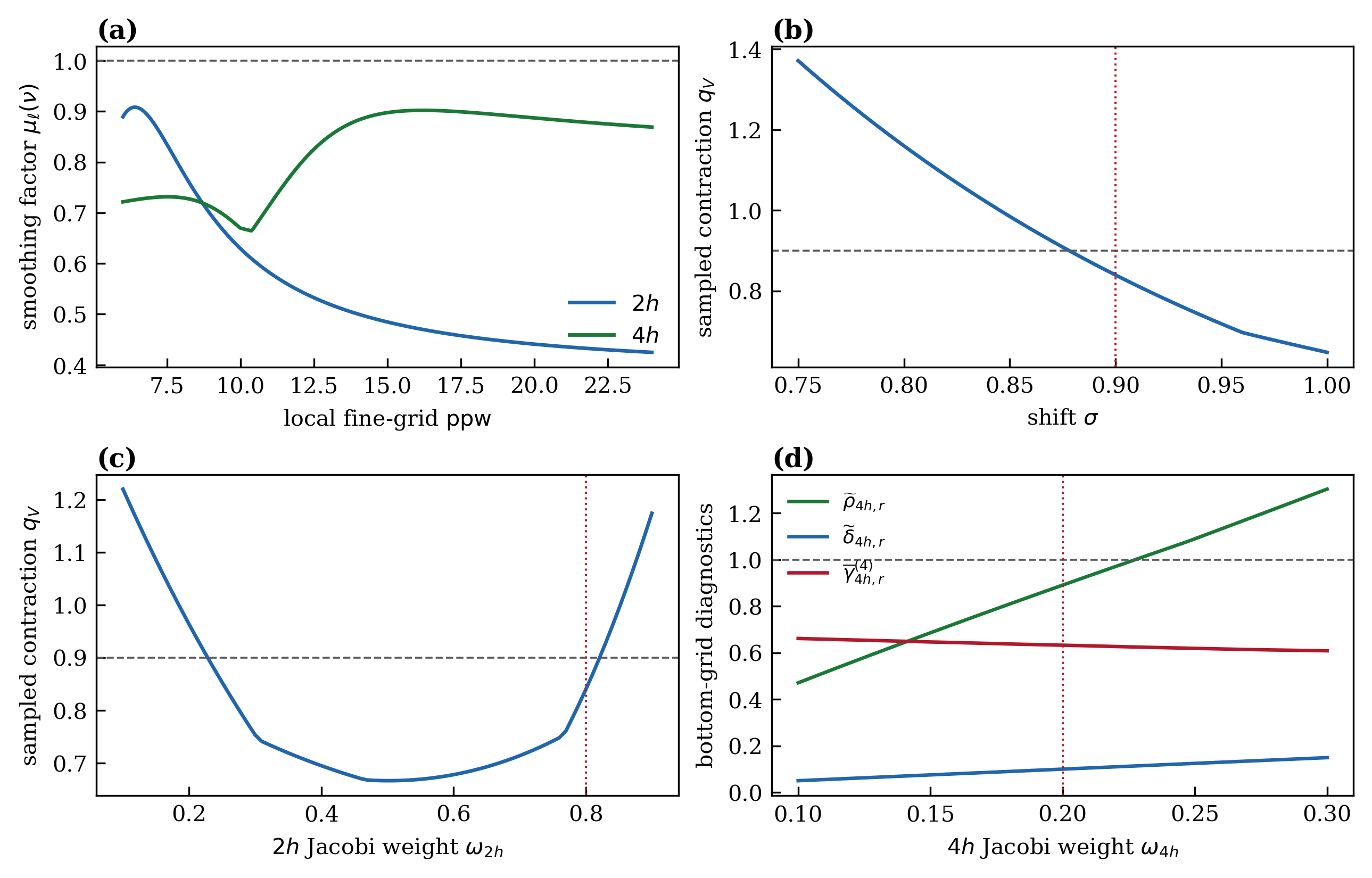}
  \caption{Source-independent Fourier diagnostics for
  \(\ppw_{\min}=6\), \(c_{\max}/c_{\min}\leq4\), and
  \(n_{J,2h}=n_{J,4h}=2\).
  (a) High-frequency factors on the shifted \(2h\) and \(4h\) levels over the
  local fine-grid sampling interval for
  \((\sigma,\omega_{2h},\omega_{4h})=(0.9,0.8,0.2)\).  (b) Numerical estimate
  \(\widetilde q_{V,r}\) versus \(\sigma\) for \(\omega_{2h}=0.8\).  (c) The
  same estimate versus \(\omega_{2h}\) for \(\sigma=0.9\).  (d) Computed
  outer radius \(\widetilde\rho_{4h,r}\), origin separation
  \(\widetilde\delta_{4h,r}\), and attained degree-four residual factor
  \(\overline{\gamma}_{4h,r}^{(4)}\) of the Jacobi-preconditioned bottom-grid
  symbol versus \(\omega_{4h}\).  Dashed lines mark reference thresholds and
  dotted lines the selected parameters.}
  \label{fig:olfd-parameter-selection}
\end{figure}

To evaluate these diagnostics numerically, let
\(\mathcal G_r\subset[6,24]\times\Theta_{\rm low}\) and
\(\mathcal H_r\subset[6,24]\times(-\pi,\pi]^3\) be the nested finite
angle--sampling sets at refinement level \(r\), and let \(\mathcal K_r\) be
the corresponding set used for the polynomial optimization.  We distinguish their
computed extrema from the continuous diagnostics by writing, with parameter
dependence suppressed,
\begin{equation}
\begin{aligned}
  \widetilde q_{V,r}
  &=\max_{(\nu,\bm\theta)\in\mathcal G_r}
    \|\widehat{\bm E}_{2h}^{\sigma}(\bm\theta;\nu)\|_2,\\
  \widetilde\rho_{4h,r}
  &=\max_{(\nu,\bm\theta)\in\mathcal H_r}
    |\widehat\lambda_{4h}^{\sigma}(\bm\theta;\nu)|,\\
  \widetilde\delta_{4h,r}
  &=\min_{(\nu,\bm\theta)\in\mathcal H_r}
    |\widehat\lambda_{4h}^{\sigma}(\bm\theta;\nu)|,\\
  \widetilde\gamma_{4h,r}^{(m)}
  &=\min_{\substack{p\in\Pi_m^{\mathbb{C}}\\p(0)=1}}
    \max_{(\nu,\bm\theta)\in\mathcal K_r}
    |p(\widehat\lambda_{4h}^{\sigma}(\bm\theta;\nu))|.
\end{aligned}
  \label{eq:sampled-parameter-diagnostics}
\end{equation}
The complex minimax problem is evaluated by a supporting-half-plane linear
program.  If \(t_r^{(m)}\) and \(p_r^{(m)}\) are its optimum and polynomial,
respectively, let \(\overline{\gamma}_{4h,r}^{(m)}\) denote the maximum of
\(|p_r^{(m)}(\widehat\lambda_{4h}^{\sigma})|\) on \(\mathcal K_r\).
For 64 uniformly spaced supporting directions,
\(t_r^{(m)}\leq\widetilde\gamma_{4h,r}^{(m)}
\leq\overline{\gamma}_{4h,r}^{(m)}
\leq\sec(\pi/64)t_r^{(m)}\), where \(\sec(\pi/64)<1.00121\).
The symmetry reduction is exact: sign changes leave all cosine symbols
unchanged, while a coordinate permutation transforms
\(\widehat{\bm E}_{2h}^{\sigma}\) by
\(\bm\Pi\widehat{\bm E}_{2h}^{\sigma}\bm\Pi^*\) and leaves
\(\widehat\lambda_{4h}^{\sigma}\) unchanged.  Thus the matrix two-norm and
scalar modulus are invariant, and one representative is retained from each
sign-change and coordinate-permutation orbit.
The same tilde convention is used below for the refined numerical estimates
of the remaining one-dimensional supremum in
\(\mu_\ell^{\rm wc}\).
The phase audit uses up to \(257\times257\) directions and 145 local-sampling
values.  After symmetry reduction, the refined \(q_V\) scan uses 28 angular
nodes per coordinate and 145 local-sampling values, while the degree-four
optimization uses 33 and 73, respectively.  Independent evaluation on 49
angular nodes and 145 local-sampling values gives a maximum modulus of 0.632,
matching \(\overline{\gamma}_{4h,r}^{(4)}\) to three digits.
The threshold crossings reported below are unchanged to two digits
under successive refinement.  Figure~\ref{fig:olfd-parameter-selection}
shows overview scans of the same diagnostics; the quoted crossings are
recomputed on the finest nested sets.  All scans use
\(n_{J,2h}=n_{J,4h}=2\), and these arguments are suppressed below.  The
one-parameter scans are cross-sections used to delineate admissible parameter
regimes.

The continuous condition \(q_V<1\) motivates using
\(\widetilde q_{V,r}\leq0.9\) as a conservative numerical admissibility
threshold for the auxiliary cycle.  At \(\omega_{2h}=0.8\), the
scan crosses this threshold near \(\sigma=0.88\); at \(\sigma=0.9\), the
corresponding interval in \(\omega_{2h}\) is approximately \(0.23\) to
\(0.82\).  These crossings delimit a sampled contractive regime for this
stationary cycle.  The two terms in
\eqref{eq:lfa-shift-balance} additionally explain why the shift should not be
increased without need.

On the bottom grid, increasing \(\omega_{4h}\) moves the computed spectrum away
from the origin but also enlarges its outer radius.  To maintain a comparable
spectral scale across the candidate weights, we adopt
\(\widetilde\rho_{4h,r}\leq1\) as a conservative normalization; its refined
unit-radius crossing occurs near \(\omega_{4h}=0.227\).  Within this range,
\(\overline{\gamma}_{4h,r}^{(4)}\) gives the residual factor attained by the
computed degree-four polynomial on the full sampled set.  At
\((\sigma,\omega_{4h})=(0.9,0.2)\), the refined values are
\(\widetilde\rho_{4h,r}=0.891\),
\(\widetilde\delta_{4h,r}=0.100\), and
\(\overline{\gamma}_{4h,r}^{(4)}=0.632\).  Thus the selected weight remains
below the radial crossing while admitting substantial residual reduction by a
degree-four polynomial.

For the full \(h\)-to-\(2h\)-to-\(4h\) hierarchy, write
\(\bm q=\bm r+2\bm p\), with
\(\bm p,\bm r\in\{0,1\}^3\).  Composing the two aliasing stages shows that the
corresponding harmonic space is spanned by exactly the 64 modes
\begin{equation}
  \left\{\bm\theta+\frac{\pi}{2}\bm q:
  \bm q\in\{0,1,2,3\}^3\right\},
  \qquad \bm\theta\in(-\pi/4,\pi/4]^3.
  \label{eq:three-grid-harmonics}
\end{equation}
The stationary stencils and smoothers act diagonally on these modes, while
restriction and prolongation couple only aliases within the same set.  If each
Krylov map is replaced by a prescribed polynomial whose coefficients are
independent of the right-hand side, the resulting linear stationary map leaves
this span invariant and reduces to \(64\times64\) Fourier blocks, recording the
two-stage aliasing pattern.  The phase scan motivates the rediscretized
physical coarse operator through its quantitative fine--coarse agreement.
Independently, the contraction
and finite-polynomial bottom diagnostics define source-independent admissibility
criteria for the auxiliary hierarchy.

\section{Numerical Validation}
\label{sec:validation}

Before turning to heterogeneous media and large-scale calculations, we use a
homogeneous point-source problem to \revision{establish} \revision{wavefield} \revision{accuracy} and \revision{choose} \revision{a}
stopping tolerance and arithmetic precision \revision{consistent} \revision{with} \revision{that} \revision{accuracy.}
The outgoing Green function provides a common \revision{reference,} \revision{so} the \revision{effect} \revision{of}
\revision{iterative} \revision{error} and \revision{rounding} \revision{can} \revision{be} \revision{assessed} \revision{against} \revision{the} \revision{accuracy} \revision{of} \revision{the}
\revision{discretized} \revision{radiating} \revision{problem.}  \revision{Guided} \revision{by} \revision{the} \revision{admissibility} diagnostics in
Section~\ref{sec:three-grid-fourier}\revision{,} \revision{we} \revision{select} \(\sigma=0.9\),
\(\omega_{2h}=0.8\), and \(\omega_{4h}=0.2\)\revision{.}  \revision{These} \revision{auxiliary} \revision{parameters}
\revision{are} kept fixed \revision{throughout} \revision{the} \revision{numerical} \revision{experiments.}

Throughout the numerical results,
mesh dimensions report the number of fine-grid intervals in each coordinate
direction, including those in the PML.  Thus, an \(n^3\) mesh has \((n+1)^3\)
nodal unknowns.  Every reported residual is evaluated for the original unshifted
system using \eqref{eq:relative-residual}.

\subsection{Green-Function Test}

The \revision{Green-function} test \revision{assesses} \revision{how} \revision{accurately} the \revision{discrete} \revision{wavefield} \revision{retains}
\revision{phase} \revision{and} \revision{amplitude} \revision{over} \revision{long} propagation \revision{distances.}  In a homogeneous medium,
the analytical \revision{reference} \revision{makes} \revision{both} \revision{effects} \revision{directly} \revision{observable,} \revision{providing} \revision{a}
\revision{basis} \revision{for} the \revision{subsequent} tolerance and \revision{precision} \revision{comparisons.}

\revision{We} set the constant wave speed to \(c=1\) and place
the source at a grid point near the center of the domain.  Its discrete
normalization is
\begin{equation}
  (f_h)_{\bm i}
  =h^{-3}\delta_{\bm i,\bm i_s},
  \qquad b_h=Q_hf_h,
  \label{eq:point-source}
\end{equation}
where \(\bm x_{\bm i_s}=\bm x_s\).  The IOFD system
\(P_hv_h=Q_hf_h\) is solved and \(u_h=Q_hv_h\) is compared directly with the
outgoing Green function
\[
  G_k(\bm{x},\bm{x}_s)
  =\frac{\exp(\ii k |\bm{x}-\bm{x}_s|)}
         {4\pi |\bm{x}-\bm{x}_s|}.
\]
Accordingly, \((\uref)_{\bm i}=G_k(\bm x_{\bm i},\bm x_s)\).
The error sequence uses
\[
  n\in\{\revision{640,1280,1672,2560,3344,6688,10032,13376}\},
  \qquad \ppw = 6, \qquad \npml = 8,
\]
where \(\npml\) denotes the number of PML layers on each side.  The wavefield
comparisons below use \revision{\(n=640\)}.
\revision{The \(1672^3\), \(3344^3\), \(6688^3\), \(10032^3\), and \(13376^3\) tests use the CPU implementation; all other
computations in this section use NVIDIA A100 40GB GPUs.}
Following \cite{TournierEtAl2022}, we report the sum of
the relative \(\ell_1\) errors in the distance-scaled real and imaginary parts,
\begin{align*}
  \mathrm{Err}
  &=
  \frac{\|W\,\operatorname{Re}(\uref-\uh)\|_1}
       {\|W\,\operatorname{Re}\uref\|_1}
  +
  \frac{\|W\,\operatorname{Im}(\uref-\uh)\|_1}
       {\|W\,\operatorname{Im}\uref\|_1},
\end{align*}
where \revision{\(\mathcal I_h=\{\bm i:\bm x_{\bm i}\in\Omega_{\rm ph},\,
r_{\bm i}=|\bm{x}_{\bm i}-\bm{x}_s|>2\pi/k,\,
\bm i\in\mathcal S_h\}\)}, the norms are restricted to
\(\mathcal I_h\), and \(W_{\bm i\bm i}=r_{\bm i}\).  \revision{The} \revision{uniform} \revision{receiver}
\revision{lattice} \revision{\(\mathcal S_h\)} \revision{has} \revision{spacing} \revision{\(4h\)} \revision{in} \revision{the} \revision{GPU} \revision{tests} \revision{and} \revision{\(h\)} \revision{in}
\revision{the} \revision{CPU} \revision{tests.}  \revision{On} the same receiver \revision{set,} \revision{we} \revision{also} \revision{measure} \revision{the} \revision{relative}
\revision{\(L^2\)} \revision{error} \revision{of} \revision{the} \revision{complex} \revision{wavefield,}
\[
  \revision{E_{L^2}=\frac{\|\uref-\uh\|_2}{\|\uref\|_2}},
\]
\revision{without} \revision{distance} \revision{weighting.}

\begin{figure}[htbp]
\centering
\includegraphics[width=\textwidth]{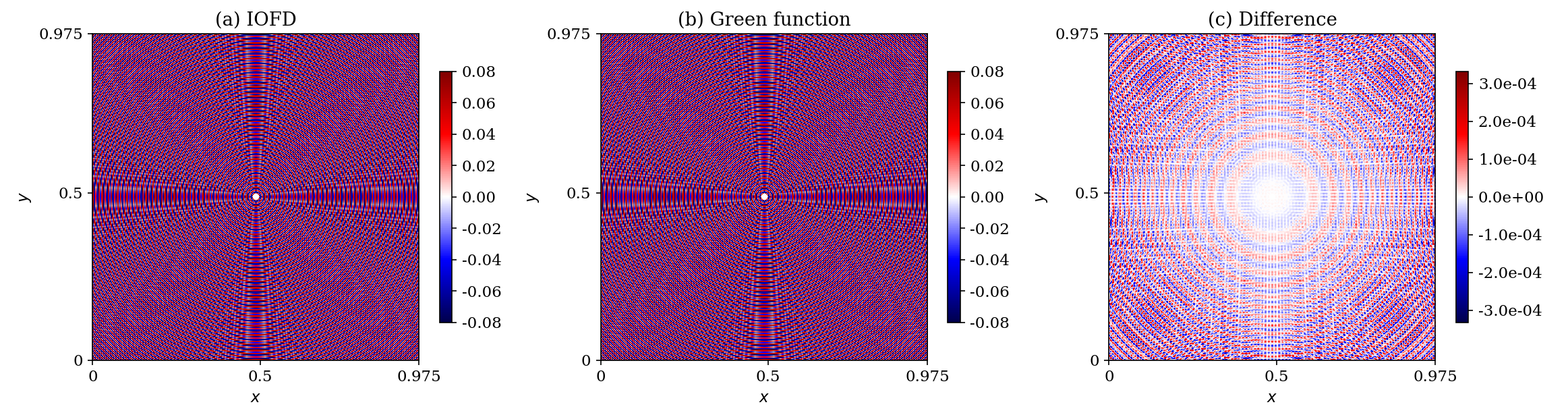}
\setlength{\abovecaptionskip}{3pt}
\caption{Homogeneous medium.  The panels show the distance-scaled IOFD
\revision{real} wavefield, the distance-scaled \revision{real} \revision{part} \revision{of} \revision{the} outgoing Green function, and their pointwise
difference on the full physical plane crossing the source.}
\label{fig:homo-green-slice}
\end{figure}

\begin{table}[htbp]
\centering
\caption{Homogeneous Green-function \revision{errors.} \revision{Mesh} \revision{dimensions} \revision{count} fine-grid intervals including the \revision{PML.}}
\label{tab:homo-green-validation}
\footnotesize

\begin{tabular}{ccc}
\toprule
\multicolumn{1}{c}{mesh} & \multicolumn{1}{c}{\(\mathrm{Err}\)} & \multicolumn{1}{c}{\revision{\(E_{L^2}\)}} \\
\midrule
\revision{\(640^3\)} & \revision{\(2.61\times10^{-3}\)} & \revision{\(1.13\times10^{-3}\)} \\
\revision{\(1280^3\)} & \revision{\(4.28\times10^{-3}\)} & \revision{\(1.78\times10^{-3}\)} \\
\revision{\(1672^3\)} & \revision{\(4.70\times10^{-3}\)} & \revision{\(1.95\times10^{-3}\)} \\
\revision{\(2560^3\)} & \revision{\(5.16\times10^{-3}\)} & \revision{\(2.27\times10^{-3}\)} \\
\revision{\(3344^3\)} & \revision{\(5.28\times10^{-3}\)} & \revision{\(2.40\times10^{-3}\)} \\
\revision{\(6688^3\)} & \revision{\(8.53\times10^{-3}\)} & \revision{\(3.74\times10^{-3}\)} \\
\revision{\(10032^3\)} & \revision{\(1.22\times10^{-2}\)} & \revision{\(5.34\times10^{-3}\)} \\
\revision{\(13376^3\)} & \revision{\(1.60\times10^{-2}\)} & \revision{\(6.99\times10^{-3}\)} \\
\bottomrule
\end{tabular}%

\end{table}

\begin{figure}[htbp]
\centering
\includegraphics[width=0.93\textwidth]{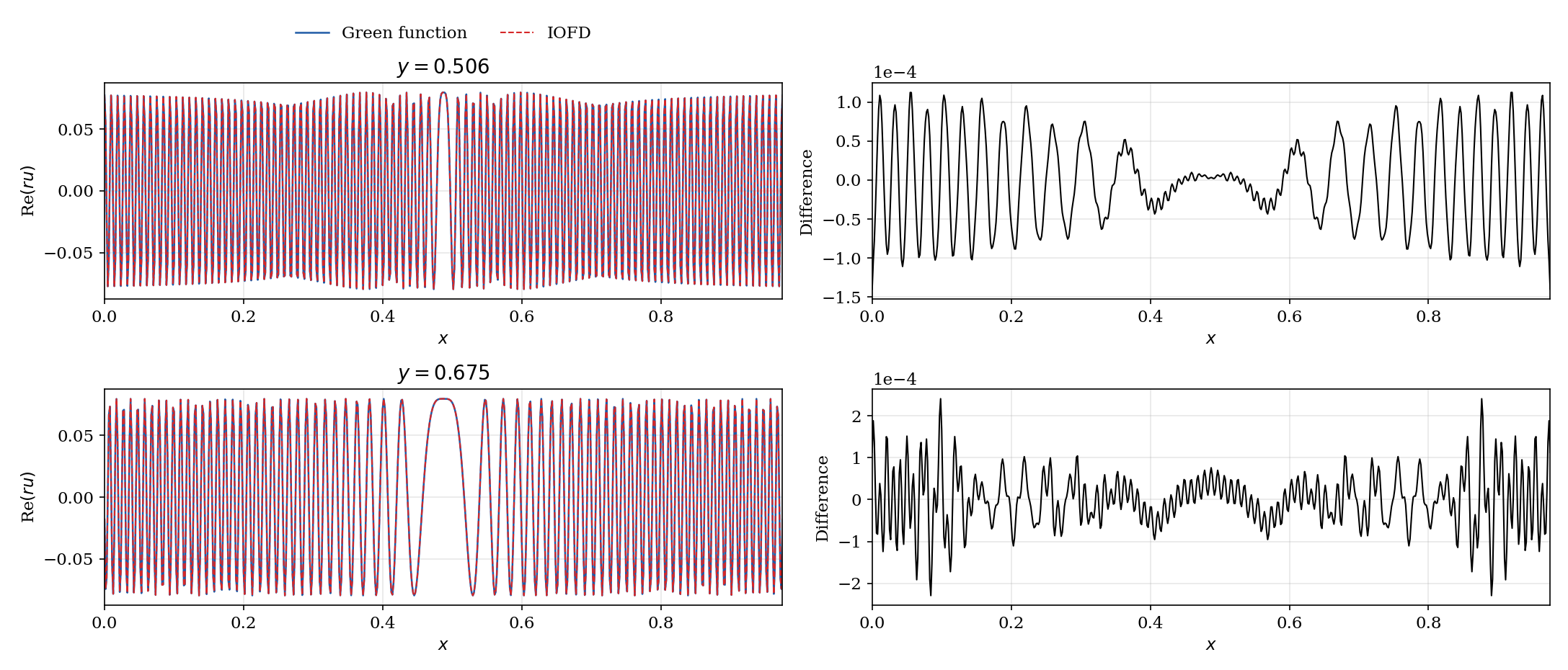}
\setlength{\abovecaptionskip}{3pt}
\caption{One-dimensional receiver-line comparisons for the homogeneous
Green-function test.  The left panels compare the distance-scaled real part of
the IOFD solution with the analytical reference.  The right panels show the
corresponding numerical difference on the same receiver \revision{lines,} \revision{at} \revision{offsets} \revision{of}
\revision{two} \revision{and} \revision{twenty} \revision{wavelengths} \revision{from} \revision{the} \revision{source.}}
\label{fig:homo-green-lines}
\end{figure}

\revision{Figure~}\revision{\ref{fig:homo-green-slice}} \revision{shows} \revision{close} \revision{agreement} \revision{between} \revision{the} numerical
and analytical \revision{wavefields} \revision{across} the physical domain.  \revision{The} receiver profiles
in Figure~\ref{fig:homo-green-lines} \revision{resolve} \revision{this} \revision{agreement} \revision{more} \revision{directly:}
\revision{both} the \revision{phase} and the \revision{distance-scaled} \revision{amplitude} \revision{are} \revision{retained} over the full
propagation \revision{interval.}  \revision{Consistently,} Table~\ref{tab:homo-green-validation}
gives \revision{relative} \revision{errors} \revision{below} \revision{\(5.2\times10^{-3}\)} \revision{for} the \revision{GPU} mesh
\revision{sequence,} \revision{as} the \revision{physical} \revision{domain} \revision{spans} \revision{104} \revision{to} \revision{424} wavelengths \revision{per}
coordinate direction.  \revision{The} \revision{CPU} \revision{tests} \revision{extend} \revision{this} \revision{comparison} \revision{to} \revision{\(13376^3\)}\revision{,}
\revision{retaining} \revision{a} \revision{relative} \revision{\(L^2\)} \revision{error} \revision{of} \revision{\(6.99\times10^{-3}\)} \revision{across}
\revision{approximately} \revision{2227} \revision{wavelengths} \revision{per} \revision{coordinate} \revision{direction.}  \revision{The} \revision{small} \revision{errors} \revision{over} \revision{these} \revision{propagation} \revision{distances}
\revision{support} the \revision{use} \revision{of} \revision{six} \revision{points} \revision{per} \revision{wavelength} \revision{in} \revision{the} solver \revision{experiments.}

\subsection{Residual Tolerance}

The same Green-reference test is used to separate Krylov error from the
\revision{remaining} \revision{wavefield} \revision{error} and \revision{identify} \revision{an} \revision{appropriate} \revision{stopping} \revision{tolerance.}
We solve the homogeneous \revision{\(640^3\)}
problem on one NVIDIA A100 40GB GPU with outer relative residual tolerances
\(10^{-2}\), \(10^{-3}\), \(10^{-4}\), and \(10^{-5}\), keeping the
discretization and three-grid preconditioner fixed.
Table~\ref{tab:homo-tolerance} reports the relative
residual, Green-function error, and the relative change of the distance-scaled
receiver profile with respect to the \(10^{-5}\) solution.
Let \(\Gamma\) denote the source-crossing receiver line outside the
one-wavelength exclusion region, and write \(u_h^{(\varepsilon)}\) for the
solution obtained with outer tolerance \(\varepsilon\).  The reported profile
change is
\begin{equation}
\begin{aligned}
  \Delta_{\rm prof}(\varepsilon)
  ={}&\frac{\|W\,\operatorname{Re}(u_h^{(\varepsilon)}-u_h^{(10^{-5})})\|_{1,\Gamma}}
           {\|W\,\operatorname{Re}u_h^{(10^{-5})}\|_{1,\Gamma}}\\
  &+\frac{\|W\,\operatorname{Im}(u_h^{(\varepsilon)}-u_h^{(10^{-5})})\|_{1,\Gamma}}
           {\|W\,\operatorname{Im}u_h^{(10^{-5})}\|_{1,\Gamma}},
\end{aligned}
\label{eq:profile-change}
\end{equation}
where \(\|\cdot\|_{1,\Gamma}\) is the discrete \(\ell_1\) norm restricted to
\(\Gamma\).
Here and below, PC calls denote the number of preconditioner applications,
solve is the elapsed solution time, and solve/PC is the solve time divided
by the number of PC calls.

\begin{table}[htbp]
\centering
\caption{Effect of the outer residual tolerance for the homogeneous Green-function \revision{test} \revision{on} \revision{a} \revision{\(640^3\)} \revision{mesh.}}
\label{tab:homo-tolerance}
\footnotesize

\begin{tabular}{cccccc}
\toprule
\multicolumn{1}{c}{outer tolerance} & \multicolumn{1}{c}{PC calls} & \multicolumn{1}{c}{\(\|A_h\vh-b_h\|_2/\|b_h\|_2\)} & \multicolumn{1}{c}{solve} & \multicolumn{1}{c}{\(\mathrm{Err}\)} & \multicolumn{1}{c}{\(\Delta_{\rm prof}\)} \\
\midrule
\(10^{-2}\) & \revision{11} & \revision{\(7.50\times10^{-3}\)} & \revision{3.714} s & \revision{\(1.74\times10^{-1}\)} & \revision{\(1.14\times10^{-1}\)} \\
\(10^{-3}\) & \revision{19} & \revision{\(6.87\times10^{-4}\)} & \revision{6.428} s & \revision{\(7.52\times10^{-3}\)} & \revision{\(5.21\times10^{-3}\)} \\
\(10^{-4}\) & \revision{24} & \revision{\(7.37\times10^{-5}\)} & \revision{8.127} s & \revision{\(2.61\times10^{-3}\)} & \revision{\(1.23\times10^{-3}\)} \\
\(10^{-5}\) & \revision{31} & \revision{\(8.20\times10^{-6}\)} & \revision{10.515} s & \revision{\(2.81\times10^{-3}\)} & \revision{--} \\
\bottomrule
\end{tabular}%

\end{table}

\begin{figure}[H]
\centering
\includegraphics[width=0.93\textwidth]{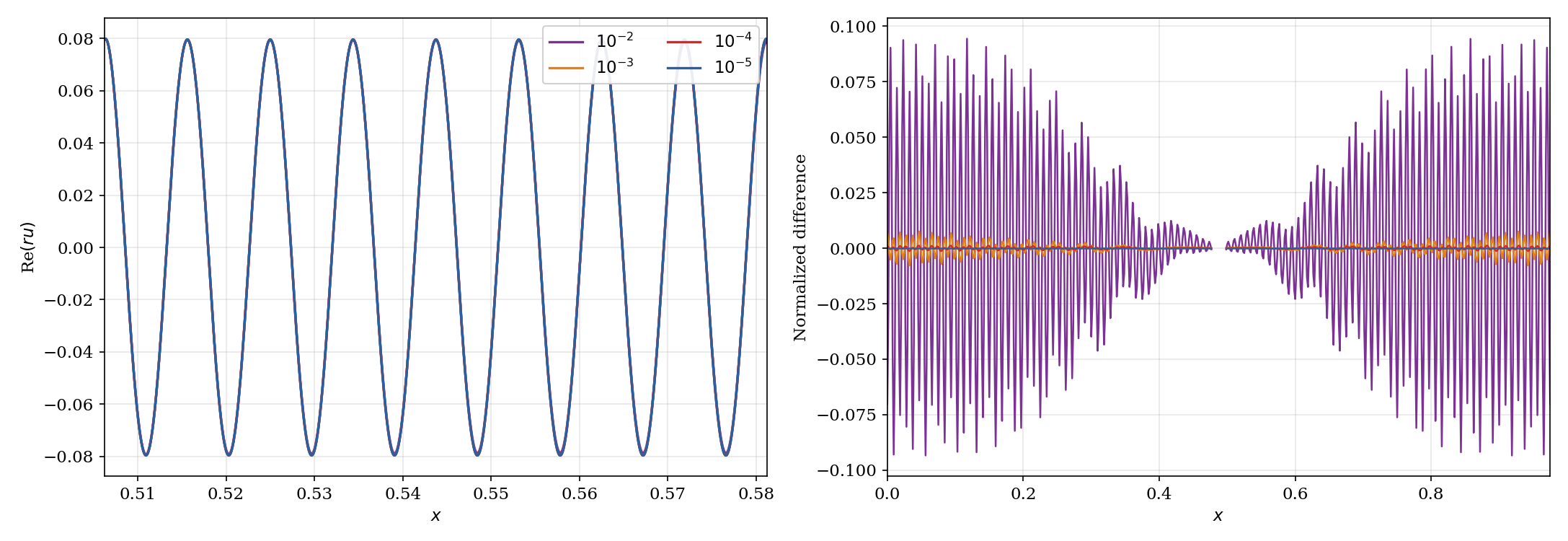}
\setlength{\abovecaptionskip}{3pt}
\caption{Effect of the outer residual tolerance on the source-crossing
receiver profile.  The left panel shows a zoom of the distance-scaled real
part of the wavefield.  The right panel shows the profile difference relative
to the \(10^{-5}\) solve, normalized by the maximum amplitude of the
\(10^{-5}\) profile.}
\label{fig:homo-tolerance-profiles}
\end{figure}

\revision{Table~}\revision{\ref{tab:homo-tolerance}} \revision{shows} \revision{that} \revision{tightening} the tolerance from
\(10^{-2}\) to \(10^{-4}\) \revision{reduces} the Green-function error by \revision{nearly} two
\revision{orders} of \revision{magnitude.}  A further reduction \revision{to} \revision{\(10^{-5}\)} \revision{requires} \revision{seven}
\revision{additional} \revision{preconditioner} \revision{applications,} \revision{while} the \revision{error} \revision{remains} \revision{at} \revision{the} \revision{same}
\revision{level.}  \revision{Figure~}\revision{\ref{fig:homo-tolerance-profiles}} \revision{confirms} \revision{that} \revision{the}
\revision{\(10^{-4}\)} \revision{and} \revision{\(10^{-5}\)} \revision{receiver} \revision{profiles} \revision{are} \revision{nearly} \revision{indistinguishable.}
\revision{We} \revision{therefore} \revision{use} \revision{\(10^{-4}\)} \revision{in} \revision{the} \revision{subsequent} \revision{experiments,} \revision{balancing}
\revision{wavefield} \revision{accuracy} \revision{against} \revision{solution} \revision{cost.}

\subsection{Arithmetic Precision}

\revision{With} \revision{the} \revision{stopping} tolerance \revision{fixed,} \revision{we} \revision{examine} \revision{whether} \revision{double} \revision{precision} \revision{yields}
a \revision{wavefield} \revision{improvement} \revision{commensurate} \revision{with} \revision{its} \revision{additional} \revision{cost.}  \revision{Complex} single
and double precision \revision{are} \revision{compared} on the same \revision{\(640^3\)} problem using two
NVIDIA A100 40GB \revision{GPUs,} \revision{with} \revision{identical} discretization, source, PML, \revision{and}
three-grid \revision{parameters.}  \revision{In} the \revision{double-precision} \revision{calculation,} stencil
coefficients, solution vectors, reductions, and Krylov \revision{orthogonalization} \revision{all}
\revision{use} \revision{double} \revision{arithmetic;} \revision{in} \revision{the} \revision{single-precision} \revision{calculation} \revision{all} \revision{vectors} \revision{and} \revision{halo} \revision{messages} \revision{are} \revision{kept} \revision{in} \revision{single} \revision{precision.}

\begin{table}[htbp]
\centering
\caption{Single- and double-precision comparison on \revision{a} \revision{\(640^3\)} \revision{mesh} \revision{using} two NVIDIA A100 40GB GPUs. Memory is the total peak over both GPUs.}
\label{tab:homo-precision}
\footnotesize

\begin{tabular}{cccccc}
\toprule
\multicolumn{1}{c}{precision} & \multicolumn{1}{c}{PC calls} & \multicolumn{1}{c}{solve} & \multicolumn{1}{c}{solve/PC} & \multicolumn{1}{c}{\revision{\(\mathrm{Err}\)}} & \multicolumn{1}{c}{memory} \\
\midrule
single & \revision{24} & \revision{4.346} s & \revision{0.181} s & \revision{\(2.61\times10^{-3}\)} & \revision{33.29} GiB \\
double & \revision{24} & \revision{8.433} s & \revision{0.351} s & \revision{\(2.61\times10^{-3}\)} & \revision{65.70} GiB \\
\bottomrule
\end{tabular}%

\end{table}

The two \revision{arithmetic} \revision{choices} require \revision{the} \revision{same} \revision{number} \revision{of} \revision{preconditioner}
\revision{applications} \revision{and} \revision{give} \revision{almost} identical \revision{Green-function} \revision{errors}
\revision{(Table~}\revision{\ref{tab:homo-precision}}\revision{).}  Double precision increases the solution
time by a factor of \revision{\(1.94\)} and \revision{nearly} \revision{doubles} the peak memory \revision{without} a
\revision{comparable} \revision{improvement} \revision{in} \revision{wavefield} \revision{accuracy.}  \revision{Single} \revision{precision} \revision{is} \revision{thus} \revision{used}
\revision{for} the remaining \revision{GPU} \revision{experiments.}

\FloatBarrier
\section{Parallel Performance}
\label{sec:scalability}

\revision{The} \revision{accuracy} \revision{tests} \revision{establish} \revision{a} \revision{common} \revision{tolerance} and \revision{precision} \revision{for} \revision{evaluating}
the cost of \revision{large-scale} \revision{wave} \revision{simulation.}  \revision{We} \revision{now} \revision{keep} \revision{the} \revision{three-grid}
\revision{configuration} \revision{fixed} \revision{across} \revision{media} \revision{with} \revision{increasingly} \revision{complex} \revision{propagation} \revision{and}
\revision{measure} \revision{its} \revision{parallel} \revision{performance} \revision{on} \revision{GPU} \revision{and} \revision{CPU} \revision{clusters.}  \revision{Weak} \revision{scaling}
\revision{examines} \revision{throughput} as the \revision{number} \revision{of} \revision{wavelengths} \revision{and} \revision{computing} \revision{resources} \revision{increase}
\revision{together,} \revision{while} strong \revision{scaling} \revision{measures} \revision{the} \revision{reduction} \revision{in} \revision{solution} \revision{time} \revision{for} \revision{a}
\revision{fixed} \revision{problem.}  The Overthrust \revision{test} \revision{then} \revision{brings} \revision{these} \revision{two} \revision{measures} \revision{together}
\revision{in} \revision{a} \revision{geophysical} \revision{setting.}

\subsection{Experimental Setup}

The fine-grid resolution is fixed at \(\ppw_{\min}=6\), and each mesh has
\(\npml=8\) PML layers on every side.  Thus the shortest wavelength is resolved
by six fine-grid points.  All runs solve the unshifted IOFD equation to a
relative residual of \(10^{-4}\), and all inner iterations execute fixed work
without convergence tests.  No parameter is changed with the mesh, medium, or
\revision{processor} count.  We denote the setup and solve times by \(T_{\rm setup}\) and
\(T_{\rm solve}\), respectively; solve/PC equals \(T_{\rm solve}\) divided by
the number of PC calls.  Memory is reported in gibibytes (GiB).
Table~\ref{tab:olfd-performance-protocol} summarizes the common configuration.

\begin{table}[!t]
\centering
\caption{Fixed discretization and three-grid parameters used in the performance
experiments.}
\label{tab:olfd-performance-protocol}
\small
\renewcommand{\arraystretch}{0.95}
\begin{tabular}{p{0.30\linewidth}p{0.60\linewidth}}
\toprule
\multicolumn{1}{c}{item} & \multicolumn{1}{c}{setting} \\
\midrule
spatial discretization & 27-point interpolated optimized finite-difference operator \\
resolution and PML & \(\ppw_{\min}=6,\quad \npml=8,\quad \gamma_{\max}=1.119058\) \\
target operator & unshifted fine-grid IOFD operator \(A_h\) \\
outer Krylov method & FGMRES with restart length \(m_{\rm out}=5\), applied to \(A_h\vh=b_h\) \\
outer stopping test & \(\|A_h\vh-b_h\|_2/\|b_h\|_2\le 10^{-4}\) \\
shifted hierarchy & independently discretized IOFD operators with shift \(\sigma=0.9\) \\
fine-grid smoother & one GMRES cycle with \(m_h=2\), right-preconditioned by \(n_{J,h}=2\) Jacobi sweeps, before and after coarse correction; \(\omega_h=0.8\) \\
\(2h\) coarse equation & \(c_{2h}=2\) FGMRES cycles with restart length \(m_{2h}=10\) on the unshifted \(A_{2h}\) equation \\
\(2h\) inner preconditioner & shifted two-grid cycle with \(n_{J,2h}=2\) Jacobi pre- and post-sweeps; \(\omega_{2h}=0.8\) \\
\(4h\) bottom correction & \(m_{4h}=4\) CA-GMRES steps, right-preconditioned by \(n_{J,4h}=2\) Jacobi sweeps; \(\omega_{4h}=0.2\) \\
precision & complex single precision \\
\bottomrule
\end{tabular}
\end{table}

The shift is confined to the auxiliary hierarchy; the outer equation,
stopping criterion, and reported residual always use the unshifted physical
operator.  All hierarchy operators and grid transfers are evaluated
\revision{matrix-free.}  \revision{The} \revision{GPU} \revision{implementation} \revision{retains} \revision{vectors} on \revision{device} \revision{and} uses a
structured three-dimensional brick decomposition, with directional
nearest-neighbor halo exchanges overlapped with interior stencil work.
\revision{The GPU results for the homogeneous, lens, wedge, and barrier models
(Tables~\ref{tab:homo-green-validation}--\ref{tab:homo-precision}, \ref{tab:olfd-a100-mesh-scaling},
\ref{tab:olfd-a100-strong-scaling-640}, and~\ref{tab:olfd-a100-strong-scaling-1280}, and
Figures~\ref{fig:homo-green-slice}--\ref{fig:homo-tolerance-profiles}; the double-precision run in
Table~\ref{tab:homo-precision} uses a double-precision build of the same code; for the lens, wedge, and barrier models the point source is at \((0.5,0.5,0.1)\) and the media are evaluated in coordinates measured from the corner of the physical box with \(h=1/n\)) and for the Overthrust model
(Table~\ref{tab:overthrust-scaling} and Figure~\ref{fig:overthrust-sections}; velocity file read with the depth
index fastest, trilinear refinement, point source at depth \(500~\mathrm{m}\) and
\((2.5,2.5)~\mathrm{km}\) laterally in the physical box) were
obtained with an independent CUDA implementation of the same discretization and
solver parameters (one MPI rank per GPU, host-staged halo exchange,
half-precision halo messages inside the preconditioner, the second-level Krylov basis and the stored outer search directions kept in block-scaled half precision with one power-of-two scale per 32 complex values, and the media
coordinates normalized over the whole box including the PML).}

\revision{The GPU experiments use nodes equipped with four NVIDIA A100 40GB
GPUs, with intra-node GPU communication routed through PCIe and the host
interconnect.} \revision{The CPU experiments use Intel Xeon Platinum 8358P nodes with 64 cores and 512 GB of memory per node.  All CPU results (Tables~\ref{tab:cpu-results}, \ref{tab:cpu-strong-scaling}, \ref{tab:overthrust-weak-scaling}, and~\ref{tab:overthrust-strong-scaling}, and the CPU rows of Table~\ref{tab:homo-green-validation}) were obtained with an independent MPI+OpenMP implementation of the same discretization and solver parameters (complex single precision, eight MPI ranks of eight OpenMP threads per node, and the outer and second-level Krylov bases and the stored outer and fine-smoother search directions kept in the same block-scaled half precision, with the outer Arnoldi coefficients formed from the exact Gram matrix of the stored basis).  GPU tables report total peak device memory, whereas peak/rank in the CPU tables denotes the maximum resident memory per MPI rank.}
\revision{The} \revision{CPU} \revision{and} \revision{GPU} \revision{implementations,} \revision{together} \revision{with} \revision{build} \revision{scripts} \revision{and}
\revision{reproduction} \revision{instructions,} \revision{are} \revision{publicly} \revision{available} \revision{at}
\revision{\url{https://github.com/aggietx/helmholtz-iofd-3g}}\revision{.}

\subsection{Benchmark Media}

Let \(\bar{\bm{x}}\in[0,1]^3\) denote coordinates normalized separately along
the sides of the physical box.  The following velocity fields are written in
these normalized coordinates.
The homogeneous medium has unit wave speed.  The smooth
converging-lens velocity field from \cite{EngquistYing2011} is
\[
 c_{\rm lens}(\bar{\bm{x}})=\frac{4}{3}\left(1-\frac{1}{2}
 \exp\{-32\|\bar{\bm{x}}-(1/2,1/2,1/2)\|_2^2\}\right).
\]
We also consider the wedge and high-contrast barrier models used in sweeping
preconditioner benchmarks \cite{PoulsonEtAl2013}.  With
\(\bar{\bm{x}}=(\bar x,\bar y,\bar z)\), their wave speeds are
\[
c_{\rm wedge}(\bar{\bm{x}})=
\begin{cases}
2,   & \bar z\leq 0.4+0.1\bar y,\\
1.5, & 0.4+0.1\bar y<\bar z\leq 0.8-0.2\bar y,\\
3,   & \bar z>0.8-0.2\bar y,
\end{cases}
\]
and
\[
c_{\rm barrier}(\bar{\bm{x}})=
\begin{cases}
10^{10}, & \bar{\bm{x}}\in[0,1]\times[0.25,0.3]\times[0,0.75],\\
1, & \text{otherwise}.
\end{cases}
\]
The calculations below retain, without modification, the auxiliary parameters
\revision{selected} \revision{in} \revision{Section~}\revision{\ref{sec:validation}}\revision{.}  \revision{This} \revision{benchmark} \revision{sequence} \revision{tests} the
\revision{same} \revision{hierarchy} \revision{under} uniform \revision{propagation,} smooth focusing, discontinuous
\revision{refraction,} and \revision{high-contrast} \revision{scattering.}

\begin{figure}[!htbp]
\centering
\includegraphics[width=\textwidth]{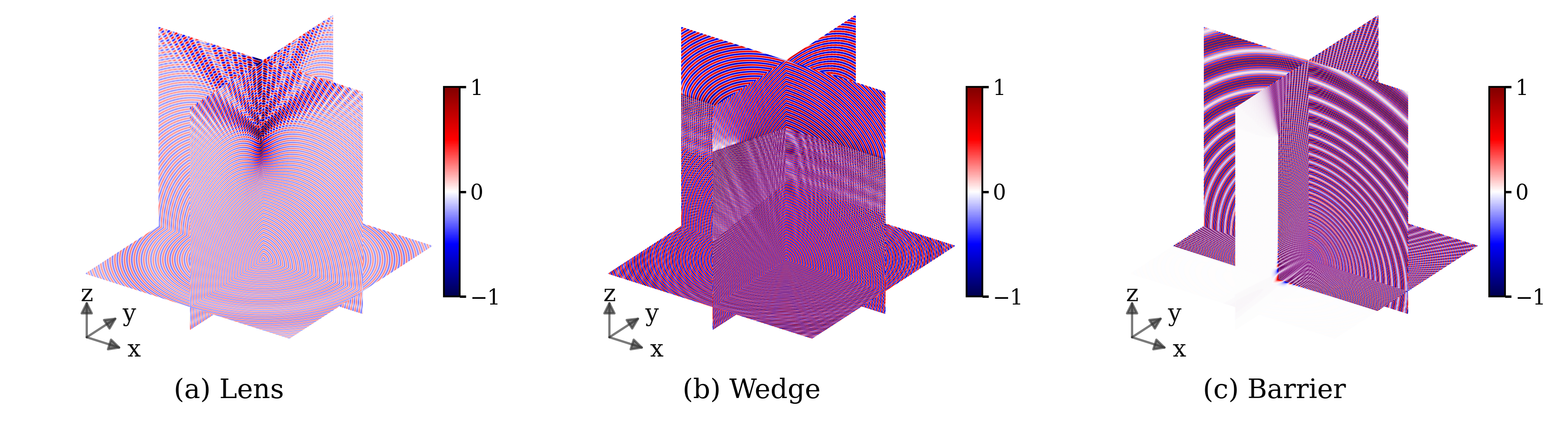}\par
\setlength{\abovecaptionskip}{2pt}
\caption{Wavefields on \revision{\(640^3\)} total meshes for the (a) converging-lens,
(b) wedge, and (c) high-contrast barrier models.  Each panel shows
\revision{\(r\operatorname{Re}(\uh)\)} on three orthogonal \revision{planes} \revision{crossing} \revision{the} \revision{source;}
PML layers are omitted.  \revision{Each} \revision{color} \revision{scale} \revision{is} \revision{normalized} \revision{by} \revision{the} \revision{99.7th}
\revision{percentile} \revision{of} \revision{the} \revision{absolute} \revision{distance-scaled} \revision{amplitude} \revision{on} \revision{the} \revision{displayed} \revision{planes.}}
\label{fig:synthetic-source-slices}
\end{figure}
\FloatBarrier

Figure~\ref{fig:synthetic-source-slices} illustrates distinct propagation
regimes under the same solver configuration.
The smooth lens bends the wavefronts into a pronounced focusing pattern, the
wedge produces refraction across oblique interfaces, and the barrier generates
strong reflection together with a downstream shadow.
Thus the benchmark set ranges from smooth phase variation to discontinuous and
high-contrast scattering before the geophysical test considered below.

\subsection{Weak Scaling}

The mesh-resource sequence uses \revision{\(640^3\)} on one GPU, \revision{\(1280^3\)} on eight
GPUs, and \revision{\(2560^3\)} on sixty-four GPUs. Let \(N_n=(n+1)^3\) \revision{and} \revision{let}
\revision{\(\tau_{n,p}\)} denote the solve time per preconditioner \revision{application.} \revision{We} \revision{use}
\begin{equation}\revision{
  \eta_{\rm w}(n,p)
  =\frac{\tau_{n_0,p_0}}{\tau_{n,p}}\frac{N_n/p}{N_{n_0}/p_0},
  \label{eq:weak-scaling-efficiency}
}\end{equation}
\revision{which} accounts for the \revision{exact} number of \revision{nodal} unknowns per \revision{processing} \revision{unit.}
\revision{Here} \revision{\((n_0,p_0)=(640,1)\)} \revision{for} \revision{GPUs} \revision{and} \revision{\((1672,64)\)} \revision{for} \revision{CPU} \revision{cores.}

\begin{table}[htbp]
\centering
\caption{\revision{Weak} scaling on NVIDIA A100 40GB \revision{GPUs,} relative to the \revision{\(640^3\)} one-GPU case for each \revision{medium.} \revision{Memory} is the total peak over all GPUs.}
\label{tab:olfd-a100-mesh-scaling}
\footnotesize
\begin{adjustbox}{max width=\textwidth}
\begin{tabular}{ccccccccc}
\toprule
\multicolumn{1}{c}{medium} & \multicolumn{1}{c}{mesh} & \multicolumn{1}{c}{PC calls} & \multicolumn{1}{c}{GPUs} & \multicolumn{1}{c}{\(T_{\rm setup}\)} & \multicolumn{1}{c}{\(T_{\rm solve}\)} & \multicolumn{1}{c}{solve/PC} & \multicolumn{1}{c}{\(\eta_{\rm w}\)} & \multicolumn{1}{c}{memory} \\
\midrule
homogeneous & \revision{\(640^3\)} & \revision{24} & 1 & \revision{0.009} s & \revision{8.120} s & \revision{0.338} s & 100.0\% & \revision{28.64} GiB \\
homogeneous & \revision{\(1280^3\)} & \revision{38} & 8 & \revision{0.028} s & \revision{12.786} s & \revision{0.337} s & \revision{100.3}\% & \revision{229.80} GiB \\
homogeneous & \revision{\(2560^3\)} & \revision{64} & 64 & \revision{0.075} s & \revision{21.684} s & \revision{0.339} s & \revision{99.5}\% & \revision{1839.39} GiB \\
\addlinespace
lens & \revision{\(640^3\)} & \revision{15} & 1 & \revision{0.009} s & \revision{5.809} s & \revision{0.387} s & 100.0\% & \revision{32.06} GiB \\
lens & \revision{\(1280^3\)} & \revision{29} & 8 & \revision{0.027} s & \revision{11.515} s & \revision{0.397} s & \revision{97.3}\% & \revision{257.05} GiB \\
lens & \revision{\(2560^3\)} & \revision{52} & 64 & \revision{0.077} s & \revision{21.217} s & \revision{0.408} s & \revision{94.6}\% & \revision{2058.11} GiB \\
\addlinespace
wedge & \revision{\(640^3\)} & \revision{36} & 1 & \revision{0.009} s & \revision{13.911} s & \revision{0.386} s & 100.0\% & \revision{32.06} GiB \\
wedge & \revision{\(1280^3\)} & \revision{69} & 8 & \revision{0.027} s & \revision{27.228} s & \revision{0.395} s & \revision{97.7}\% & \revision{257.05} GiB \\
wedge & \revision{\(2560^3\)} & \revision{139} & 64 & \revision{0.079} s & \revision{56.311} s & \revision{0.405} s & \revision{95.0}\% & \revision{2058.10} GiB \\
\addlinespace
barrier & \revision{\(640^3\)} & \revision{36} & 1 & \revision{0.009} s & \revision{13.908} s & \revision{0.386} s & 100.0\% & \revision{32.06} GiB \\
barrier & \revision{\(1280^3\)} & \revision{60} & 8 & \revision{0.027} s & \revision{23.672} s & \revision{0.395} s & \revision{97.7}\% & \revision{257.05} GiB \\
barrier & \revision{\(2560^3\)} & \revision{107} & 64 & \revision{0.074} s & \revision{43.343} s & \revision{0.405} s & \revision{95.0}\% & \revision{2058.11} GiB \\
\bottomrule
\end{tabular}%
\end{adjustbox}
\end{table}

\revision{Table~}\revision{\ref{tab:olfd-a100-mesh-scaling}} \revision{shows} \revision{that} the \revision{average} \revision{solve} \revision{time} per
\revision{preconditioner} \revision{application} remains nearly constant \revision{from} \revision{one} \revision{to} \revision{eight} \revision{GPUs,}
\revision{with} \revision{efficiencies} \revision{of} \revision{97.3--100.3}\revision{\%} \revision{across} \revision{the} \revision{four} \revision{media.}  \revision{Increasing} \revision{the}
\revision{mesh} \revision{volume} \revision{and} \revision{GPU} \revision{count} \revision{by} \revision{another} \revision{factor} \revision{of} \revision{eight} \revision{retains}
\revision{94.6--99.5}\revision{\%} \revision{efficiency} \revision{relative} \revision{to} \revision{the} \revision{one-GPU} \revision{baseline.}  \revision{Thus} \revision{the} \revision{hierarchy}
\revision{maintains} \revision{high} \revision{throughput} as both the \revision{propagation} \revision{distance} \revision{and} \revision{distributed}
\revision{problem} \revision{size} \revision{grow,} \revision{including} \revision{in} \revision{the} \revision{discontinuous} \revision{wedge} \revision{and} \revision{high-contrast}
\revision{barrier} \revision{models.}

\revision{This} \revision{throughput} \revision{translates} \revision{into} \revision{short} \revision{solution} \revision{times} \revision{at} \revision{the} \revision{largest} \revision{scale:}
\revision{the} \revision{\(2560^3\)} \revision{homogeneous} \revision{and} \revision{lens} \revision{problems} \revision{are} \revision{solved} \revision{in} \revision{21.7} \revision{and} \revision{21.2}
\revision{seconds,} \revision{respectively,} \revision{while} \revision{the} \revision{wedge} \revision{and} \revision{barrier} \revision{problems} \revision{require} \revision{56.3} \revision{and}
\revision{43.3} \revision{seconds} \revision{with} \revision{the} \revision{same} \revision{parameters.}  \revision{Setup} \revision{remains} \revision{below} \revision{0.1} \revision{seconds} \revision{in}
\revision{every} \revision{case.}  \revision{Total} \revision{peak} \revision{memory} \revision{grows} \revision{approximately} \revision{in} \revision{proportion} \revision{to} \revision{the}
\revision{number} \revision{of} \revision{unknowns,} \revision{allowing} \revision{the} \revision{complete} \revision{three-grid} \revision{hierarchy} \revision{to} \revision{fit} \revision{on} \revision{one}
\revision{40GB} \revision{GPU} \revision{at} \revision{\(640^3\)} \revision{and} \revision{to} \revision{retain} \revision{essentially} \revision{the} \revision{same} \revision{storage} \revision{per} \revision{GPU} \revision{along}
\revision{the} \revision{mesh--resource} \revision{sequence.}

\revision{The CPU results in Table~\ref{tab:cpu-results} extend the homogeneous weak-scaling sequence to \(13376^3\). Increasing the core count from 64 to 32768 raises solve/PC only from 21.438 s to 22.555 s, retaining 94.9\% efficiency after normalization by the nodal unknown count. On 32768 cores, the \(13376^3\) system, with \(2.39\times10^{12}\) unknowns, is solved in 98.9 minutes, with a maximum resident memory of 49.84 GiB per MPI rank (eight ranks of eight threads per node).}

\begin{table}[H]
\centering
\caption{\revision{CPU} \revision{weak} \revision{scaling} \revision{for} \revision{the} \revision{homogeneous} \revision{model,} \revision{relative} \revision{to} \revision{the} \revision{\(1672^3\)}\revision{,} \revision{64-core} \revision{case.} \revision{Each} \revision{node} \revision{uses} \revision{64} \revision{CPU} \revision{cores} \revision{(8 MPI ranks of 8 OpenMP threads).}}
\label{tab:cpu-results}
\footnotesize
\begin{adjustbox}{max width=\textwidth}
\begin{tabular}{cccccccc}
\toprule
\multicolumn{1}{c}{mesh} & \multicolumn{1}{c}{\revision{CPU} \revision{cores}} & \multicolumn{1}{c}{\revision{PC} \revision{calls}} & \multicolumn{1}{c}{\revision{setup}} & \multicolumn{1}{c}{\revision{solve}} & \multicolumn{1}{c}{\revision{solve/PC}} & \multicolumn{1}{c}{\revision{\(\eta_{\rm w}\)}} & \multicolumn{1}{c}{\revision{peak/rank}} \\
\midrule
\revision{\(1672^3\)} & \revision{64} & \revision{47} & \revision{1.197} s & \revision{1007.563} s & \revision{21.438} s & \revision{100.0}\% & \revision{49.845} GiB \\
\revision{\(3344^3\)} & \revision{512} & \revision{77} & \revision{1.222} s & \revision{1680.582} s & \revision{21.826} s & \revision{98.1}\% & \revision{49.877} GiB \\
\revision{\(6688^3\)} & \revision{4096} & \revision{139} & \revision{1.214} s & \revision{3083.816} s & \revision{22.186} s & \revision{96.5}\% & \revision{49.912} GiB \\
\revision{\(10032^3\)} & \revision{13824} & \revision{202} & \revision{1.244} s & \revision{4537.461} s & \revision{22.463} s & \revision{95.3}\% & \revision{49.920} GiB \\
\revision{\(13376^3\)} & \revision{32768} & \revision{263} & \revision{1.249} s & \revision{5931.822} s & \revision{22.555} s & \revision{94.9}\% & \revision{49.836} GiB \\
\bottomrule
\end{tabular}%
\end{adjustbox}
\end{table}

\FloatBarrier

\subsection{Strong Scaling}

For a fixed \revision{mesh,} \revision{the} \revision{per-application} \revision{strong-scaling} \revision{efficiency} \revision{is}
\begin{equation}\revision{
  \eta_{\rm s}(p)=\frac{p_0\tau_{p_0}}{p\tau_p},
  \label{eq:strong-scaling-efficiency}
}\end{equation}
\revision{where} \revision{\(p_0\)} \revision{is} \revision{the} \revision{smallest} \revision{processor} \revision{count} \revision{tested} \revision{for} \revision{each} mesh \revision{on} \revision{each}
\revision{platform.} \revision{For} \revision{the} \revision{GPU} \revision{experiments,} \revision{\(p_0=1\)} \revision{for} \revision{\(640^3\)} \revision{and}
\revision{\(p_0=8\)} \revision{for} \revision{\(1280^3\)}.

\begin{table}[htbp]
\centering
\caption{Strong scaling on a fixed \revision{\(640^3\)} total mesh using \revision{NVIDIA} A100 40GB GPUs. Memory is the total peak over all GPUs.}
\label{tab:olfd-a100-strong-scaling-640}
\footnotesize
\begin{adjustbox}{max width=\textwidth}
\begin{tabular}{cccccccc}
\toprule
\multicolumn{1}{c}{medium} & \multicolumn{1}{c}{PC calls} & \multicolumn{1}{c}{GPUs} & \multicolumn{1}{c}{\(T_{\rm setup}\)} & \multicolumn{1}{c}{\(T_{\rm solve}\)} & \multicolumn{1}{c}{solve/PC} & \multicolumn{1}{c}{\(\eta_{\rm s}\)} & \multicolumn{1}{c}{memory} \\
\midrule
homogeneous & \revision{24} & 1 & \revision{0.009} s & \revision{8.120} s & \revision{0.338} s & 100.0\% & \revision{28.64} GiB \\
homogeneous & \revision{24} & 2 & \revision{0.004} s & \revision{4.253} s & \revision{0.177} s & \revision{95.5}\% & \revision{29.34} GiB \\
homogeneous & \revision{24} & 4 & \revision{0.003} s & \revision{2.402} s & \revision{0.100} s & \revision{84.5}\% & \revision{30.46} GiB \\
\addlinespace
lens & \revision{15} & 1 & \revision{0.009} s & \revision{5.809} s & \revision{0.387} s & 100.0\% & \revision{32.06} GiB \\
lens & \revision{15} & 2 & \revision{0.004} s & \revision{3.013} s & \revision{0.201} s & \revision{96.4}\% & \revision{32.78} GiB \\
lens & \revision{15} & 4 & \revision{0.004} s & \revision{1.705} s & \revision{0.114} s & \revision{85.2}\% & \revision{33.96} GiB \\
\addlinespace
wedge & \revision{36} & 1 & \revision{0.009} s & \revision{13.911} s & \revision{0.386} s & 100.0\% & \revision{32.06} GiB \\
wedge & \revision{36} & 2 & \revision{0.004} s & \revision{7.215} s & \revision{0.200} s & \revision{96.4}\% & \revision{32.78} GiB \\
wedge & \revision{36} & 4 & \revision{0.003} s & \revision{3.925} s & \revision{0.109} s & \revision{88.6}\% & \revision{33.96} GiB \\
\addlinespace
barrier & \revision{36} & 1 & \revision{0.009} s & \revision{13.908} s & \revision{0.386} s & 100.0\% & \revision{32.06} GiB \\
barrier & \revision{36} & 2 & \revision{0.004} s & \revision{7.213} s & \revision{0.200} s & \revision{96.4}\% & \revision{32.78} GiB \\
barrier & \revision{36} & 4 & \revision{0.004} s & \revision{3.970} s & \revision{0.110} s & \revision{87.6}\% & \revision{33.96} GiB \\
\bottomrule
\end{tabular}%
\end{adjustbox}
\end{table}

\begin{table}[htbp]
\centering
\caption{Strong scaling on a fixed \revision{\(1280^3\)} total mesh using \revision{NVIDIA} A100 40GB GPUs. Memory is the total peak over all GPUs.}
\label{tab:olfd-a100-strong-scaling-1280}
\footnotesize
\begin{adjustbox}{max width=\textwidth}
\begin{tabular}{cccccccc}
\toprule
\multicolumn{1}{c}{medium} & \multicolumn{1}{c}{PC calls} & \multicolumn{1}{c}{GPUs} & \multicolumn{1}{c}{\(T_{\rm setup}\)} & \multicolumn{1}{c}{\(T_{\rm solve}\)} & \multicolumn{1}{c}{solve/PC} & \multicolumn{1}{c}{\(\eta_{\rm s}\)} & \multicolumn{1}{c}{memory} \\
\midrule
homogeneous & \revision{38} & 8 & \revision{0.028} s & \revision{12.786} s & \revision{0.337} s & 100.0\% & \revision{229.80} GiB \\
homogeneous & \revision{38} & 16 & \revision{0.032} s & \revision{6.774} s & \revision{0.178} s & \revision{94.4}\% & \revision{234.94} GiB \\
homogeneous & \revision{38} & 32 & \revision{0.041} s & \revision{3.972} s & \revision{0.104} s & \revision{80.5}\% & \revision{243.47} GiB \\
\addlinespace
lens & \revision{29} & 8 & \revision{0.027} s & \revision{11.515} s & \revision{0.397} s & 100.0\% & \revision{257.05} GiB \\
lens & \revision{29} & 16 & \revision{0.032} s & \revision{6.071} s & \revision{0.209} s & \revision{94.8}\% & \revision{262.52} GiB \\
lens & \revision{29} & 32 & \revision{0.040} s & \revision{3.425} s & \revision{0.118} s & \revision{84.1}\% & \revision{271.38} GiB \\
\addlinespace
wedge & \revision{69} & 8 & \revision{0.027} s & \revision{27.228} s & \revision{0.395} s & 100.0\% & \revision{257.05} GiB \\
wedge & \revision{69} & 16 & \revision{0.032} s & \revision{14.277} s & \revision{0.207} s & \revision{95.4}\% & \revision{262.52} GiB \\
wedge & \revision{69} & 32 & \revision{0.040} s & \revision{8.009} s & \revision{0.116} s & \revision{85.0}\% & \revision{271.38} GiB \\
\addlinespace
barrier & \revision{60} & 8 & \revision{0.027} s & \revision{23.672} s & \revision{0.395} s & 100.0\% & \revision{257.05} GiB \\
barrier & \revision{60} & 16 & \revision{0.032} s & \revision{12.429} s & \revision{0.207} s & \revision{95.2}\% & \revision{262.52} GiB \\
barrier & \revision{60} & 32 & \revision{0.041} s & \revision{6.970} s & \revision{0.116} s & \revision{84.9}\% & \revision{271.38} GiB \\
\bottomrule
\end{tabular}%
\end{adjustbox}
\end{table}

\revision{Tables~}\revision{\ref{tab:olfd-a100-strong-scaling-640}}
\revision{and~}\revision{\ref{tab:olfd-a100-strong-scaling-1280}} \revision{show} \revision{consistent} \revision{strong} \revision{scaling}
\revision{at} \revision{both} \revision{baseline} \revision{sizes.}  \revision{On} \revision{the} \revision{\(640^3\)} \revision{mesh,} \revision{two} \revision{GPUs} \revision{retain}
\revision{95.5--96.4}\revision{\%} \revision{efficiency} \revision{and} \revision{four} \revision{GPUs} \revision{retain} \revision{84.5--88.6}\revision{\%}\revision{.}  \revision{On} \revision{the}
\revision{\(1280^3\)} \revision{mesh,} \revision{increasing} \revision{the} \revision{allocation} \revision{from} \revision{eight} \revision{to} \revision{sixteen} \revision{GPUs} \revision{gives}
\revision{94.4--95.4}\revision{\%} \revision{efficiency,} \revision{while} \revision{thirty-two} \revision{GPUs} \revision{retain} \revision{80.5--85.0}\revision{\%}\revision{.}
\revision{The} \revision{corresponding} \revision{reduction} \revision{in} \revision{solution} \revision{time} \revision{from} \revision{eight} \revision{to} \revision{thirty-two} \revision{GPUs}
\revision{is} \revision{a} \revision{factor} \revision{of} \revision{3.22--3.40} \revision{across} \revision{the} \revision{four} \revision{media.}  \revision{The} \revision{number} \revision{of} \revision{preconditioner}
\revision{applications} \revision{is} \revision{unchanged} \revision{within} \revision{each} \revision{fixed-mesh} \revision{test} \revision{so} \revision{these} \revision{reductions}
\revision{reflect} \revision{parallel} \revision{acceleration} \revision{of} \revision{the} \revision{same} \revision{solver} \revision{rather} \revision{than} \revision{a} \revision{change} \revision{in}
\revision{iteration} \revision{strength.}  \revision{Together} \revision{with} \revision{the} \revision{weak-scaling} \revision{results,} \revision{they} \revision{show} \revision{that}
\revision{the} \revision{hierarchy} \revision{supports} \revision{both} \revision{larger} \revision{wave} \revision{problems} \revision{and} \revision{faster} \revision{solution} \revision{of} \revision{a}
\revision{fixed} \revision{problem.}
\FloatBarrier
The CPU experiments also exhibit effective strong scaling (Table~\ref{tab:cpu-strong-scaling}). \revision{Relative to the smallest allocation for each mesh, efficiencies remain at 86.4\% for \(1672^3\) on 1024 cores, 84.3\% for \(3344^3\) on 8192 cores, and 86.7\% for \(6688^3\) on 32768 cores, where the \(6688^3\) system is solved in 7.4 minutes.  Each mesh requires the same number of preconditioner applications on every allocation (47, 77, and 139).}

\begin{table}[H]
\centering
\caption{\revision{CPU} \revision{strong} scaling \revision{for} the \revision{homogeneous} \revision{model.}}
\label{tab:cpu-strong-scaling}
\footnotesize
\begin{adjustbox}{max width=\textwidth}
\begin{tabular}{cccccccc}
\toprule
\multicolumn{1}{c}{\revision{mesh}} & \multicolumn{1}{c}{\revision{CPU} \revision{cores}} & \multicolumn{1}{c}{PC calls} & \multicolumn{1}{c}{\revision{setup}} & \multicolumn{1}{c}{\revision{solve}} & \multicolumn{1}{c}{solve/PC} & \multicolumn{1}{c}{\(\eta_{\rm s}\)} & \multicolumn{1}{c}{\revision{peak/rank}} \\
\midrule
\revision{\(1672^3\)} & \revision{64} & \revision{47} & \revision{1.197} s & \revision{1007.563} s & \revision{21.438} s & \revision{100.0}\% & \revision{49.845} GiB \\
\revision{\(1672^3\)} & \revision{128} & \revision{47} & \revision{0.611} s & \revision{518.199} s & \revision{11.026} s & \revision{97.2}\% & \revision{25.218} GiB \\
\revision{\(1672^3\)} & \revision{256} & \revision{47} & \revision{0.309} s & \revision{261.435} s & \revision{5.562} s & \revision{96.4}\% & \revision{12.790} GiB \\
\revision{\(1672^3\)} & \revision{512} & \revision{47} & \revision{0.163} s & \revision{139.417} s & \revision{2.966} s & \revision{90.3}\% & \revision{6.523} GiB \\
\revision{\(1672^3\)} & \revision{1024} & \revision{47} & \revision{0.089} s & \revision{72.889} s & \revision{1.551} s & \revision{86.4}\% & \revision{3.410} GiB \\
\addlinespace
\revision{\(3344^3\)} & \revision{512} & \revision{77} & \revision{1.222} s & \revision{1680.582} s & \revision{21.826} s & \revision{100.0}\% & \revision{49.877} GiB \\
\revision{\(3344^3\)} & \revision{1024} & \revision{77} & \revision{0.637} s & \revision{853.059} s & \revision{11.079} s & \revision{98.5}\% & \revision{25.230} GiB \\
\revision{\(3344^3\)} & \revision{2048} & \revision{77} & \revision{0.325} s & \revision{438.380} s & \revision{5.693} s & \revision{95.8}\% & \revision{12.806} GiB \\
\revision{\(3344^3\)} & \revision{4096} & \revision{77} & \revision{0.172} s & \revision{237.522} s & \revision{3.085} s & \revision{88.4}\% & \revision{6.554} GiB \\
\revision{\(3344^3\)} & \revision{8192} & \revision{77} & \revision{0.097} s & \revision{124.658} s & \revision{1.619} s & \revision{84.3}\% & \revision{3.431} GiB \\
\addlinespace
\revision{\(6688^3\)} & \revision{4096} & \revision{139} & \revision{1.214} s & \revision{3083.816} s & \revision{22.186} s & \revision{100.0}\% & \revision{49.912} GiB \\
\revision{\(6688^3\)} & \revision{8192} & \revision{139} & \revision{0.641} s & \revision{1573.614} s & \revision{11.321} s & \revision{98.0}\% & \revision{25.258} GiB \\
\revision{\(6688^3\)} & \revision{16384} & \revision{139} & \revision{0.327} s & \revision{816.670} s & \revision{5.875} s & \revision{94.4}\% & \revision{12.836} GiB \\
\revision{\(6688^3\)} & \revision{32768} & \revision{139} & \revision{0.171} s & \revision{444.410} s & \revision{3.197} s & \revision{86.7}\% & \revision{6.477} GiB \\
\bottomrule
\end{tabular}%
\end{adjustbox}
\end{table}

\FloatBarrier

\subsection{Overthrust Model}

We finally consider the SEG/EAGE Overthrust velocity model
\cite{AminzadehBracKunz1997,PoulsonEtAl2013}, a strongly heterogeneous model
containing discontinuous sedimentary layers and an overthrust fault.  The
\revision{\(20~\mathrm{km}\times20~\mathrm{km}\times4.6~\mathrm{km}\)} \revision{physical}
\revision{domain} \revision{is} \revision{sampled} \revision{at} \revision{a} \revision{spacing} \revision{of} \revision{\(25~\mathrm{m}\)}\revision{,} \revision{yielding}
\revision{\(801\times801\times185\)} \revision{nodal} velocity \revision{values.}  The wave speed ranges
from approximately \(2.179\) to \(6.000~\mathrm{km/s}\).  \revision{Refinement} by
trilinear \revision{interpolation} \revision{gives} \(1601\times1601\times369\) and
\(3201\times3201\times737\) \revision{velocity} samples.  After adding eight PML layers on each
side, the three \revision{GPU} meshes contain \(816\times816\times200\),
\(1616\times1616\times384\), and \(3216\times3216\times752\) fine-grid
intervals, respectively.

\begin{figure}[htbp]
\centering
\begin{minipage}[t]{0.49\textwidth}
  \centering
  \includegraphics[width=\linewidth]{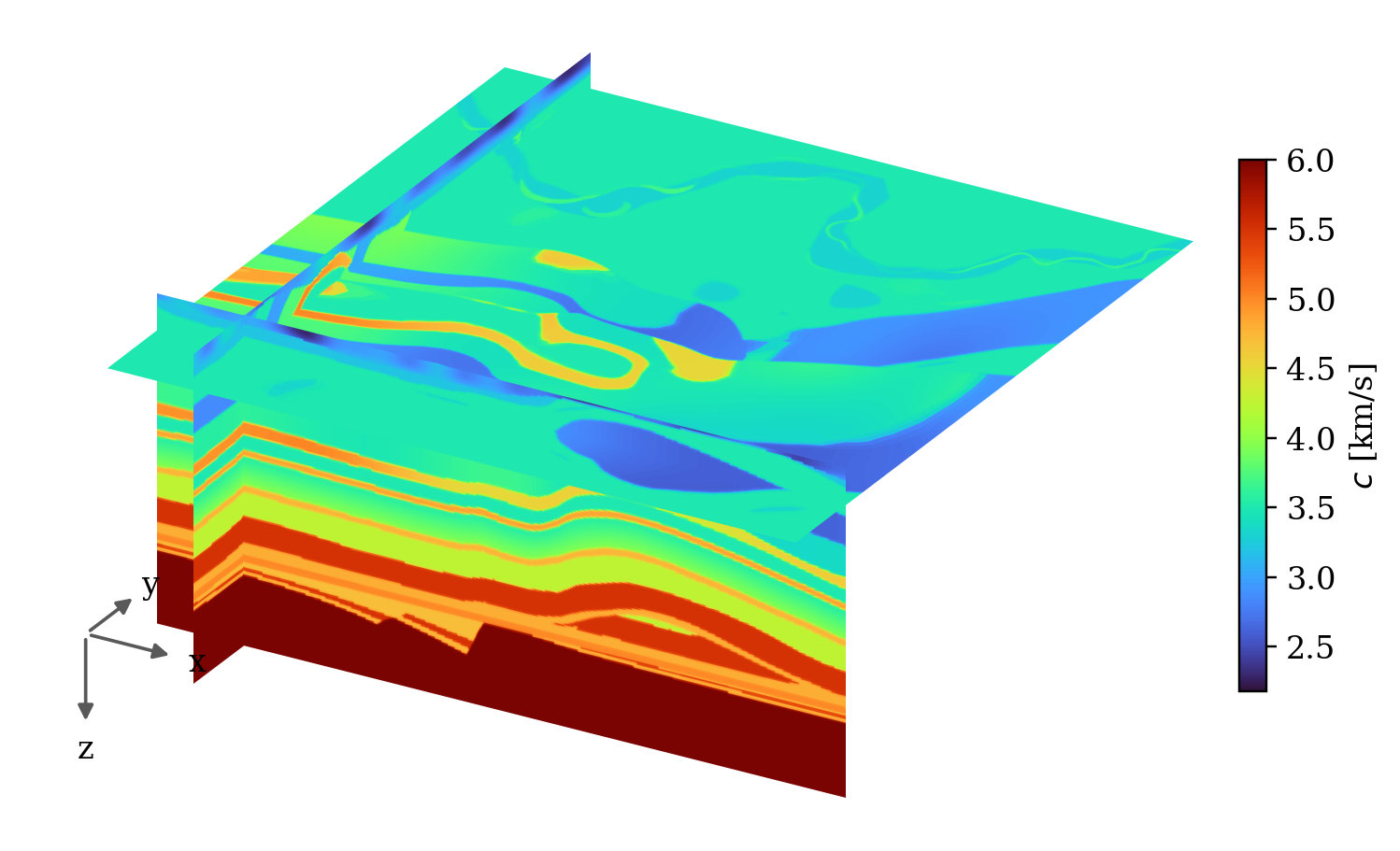}
  \par\smallskip\textup{(a)} Velocity
\end{minipage}\hfill
\begin{minipage}[t]{0.49\textwidth}
  \centering
  \includegraphics[width=\linewidth]{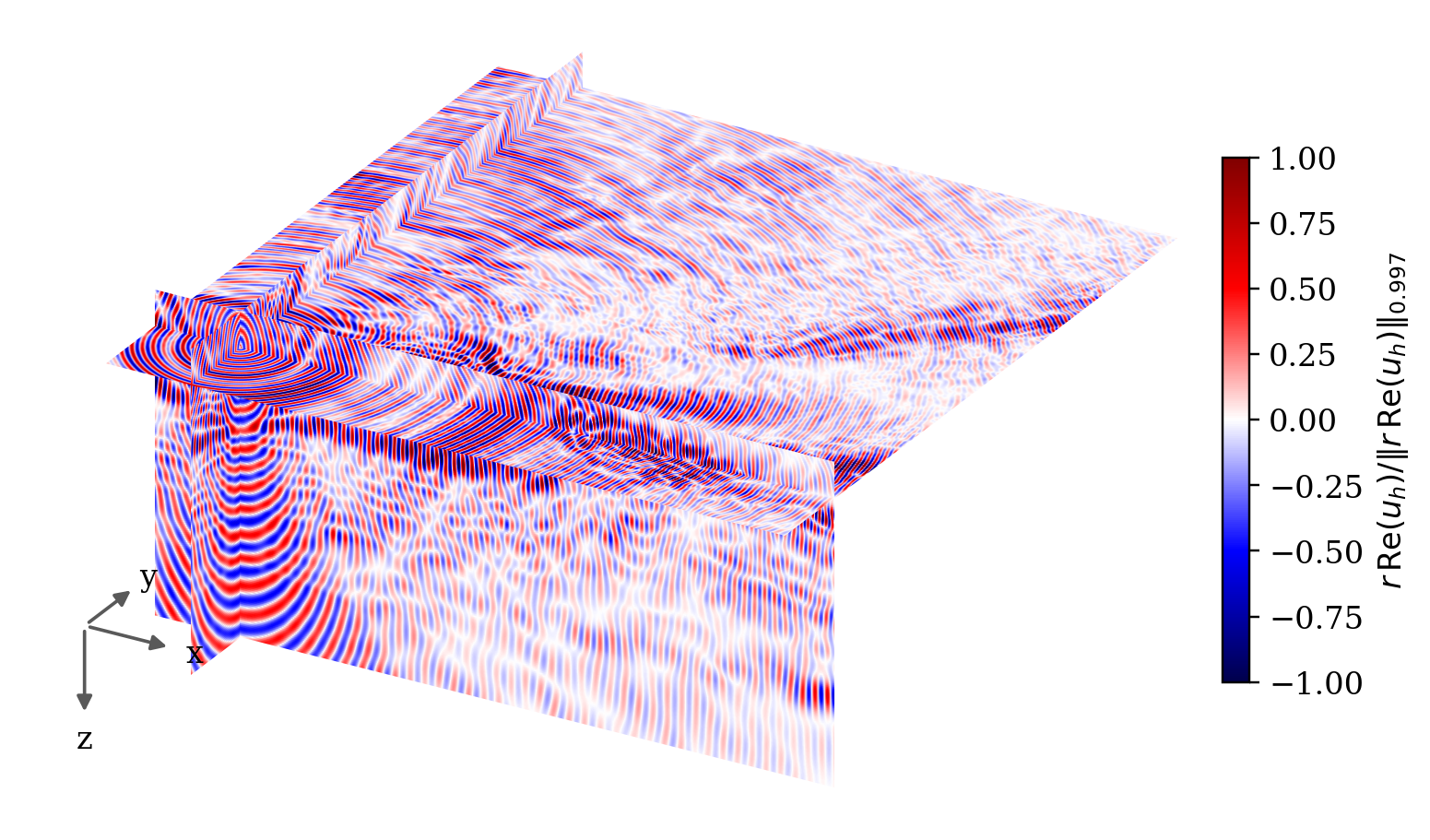}
  \par\smallskip\textup{(b)} Wavefield
\end{minipage}
\vspace{-0.4em}
\caption{Native SEG/EAGE Overthrust model: (a) velocity and (b)
distance-compensated real part of the IOFD wavefield at \(14.5~\mathrm{Hz}\)
on the same three orthogonal sections.  The wavefield color scale is normalized
by the 99.7th percentile of the distance-compensated amplitude.}
\label{fig:overthrust-sections}
\end{figure}

The velocity sections in Figure~\ref{fig:overthrust-sections}(a) expose the
layered structure, lateral variation, and fault geometry of the model.  The
corresponding wavefield in Figure~\ref{fig:overthrust-sections}(b) shows the
resulting refraction and reflection over the full survey-scale domain.  This
example therefore tests the same three-grid hierarchy when the local wavelength
and principal propagation paths vary strongly in all three dimensions.

The PML width remains fixed at eight fine-grid layers while the mesh
spacing is reduced from \(25\) to \(12.5\) and \(6.25~\mathrm{m}\).  Keeping
\(\ppw_{\min}=6\) gives frequencies \(\omega/(2\pi)\) of \(14.5\),
\(29.1\), and \(58.1~\mathrm{Hz}\), respectively.
Table~\ref{tab:overthrust-scaling}
reports \revision{convergence} \revision{and} \revision{timing} \revision{for} \revision{the} \revision{three} \revision{frequencies.}  For each fixed mesh, the strong-scaling
efficiency \(\eta_{\rm s}\) is based on solve time per preconditioner application
and the smallest GPU allocation reported for that mesh.  Memory is the total
peak over all GPUs.  Since the PML width is fixed in grid layers, successive
refinements do not produce an exact eightfold increase in total unknowns.  The
weak-scaling comparison therefore uses the \revision{same} \revision{unknown-count} normalization \revision{as}
\eqref{eq:weak-scaling-efficiency}, with
\(N=(n_x+1)(n_y+1)(n_z+1)\) for a rectangular \revision{mesh} \revision{and} \revision{the} \revision{one-GPU} \revision{native}
\revision{mesh} \revision{as} \revision{the} \revision{baseline.}

\begin{table}[htbp]
\centering
\caption{Parallel performance for the SEG/EAGE Overthrust model on NVIDIA A100 40GB GPUs. Mesh dimensions count fine-grid intervals \revision{including} \revision{the} \revision{PML.}}
\label{tab:overthrust-scaling}
\footnotesize
\begin{adjustbox}{max width=\textwidth}
\begin{tabular}{ccccccccc}
\toprule
\multicolumn{1}{c}{frequency} & \multicolumn{1}{c}{mesh} & \multicolumn{1}{c}{PC calls} & \multicolumn{1}{c}{GPUs} & \multicolumn{1}{c}{\(T_{\rm setup}\)} & \multicolumn{1}{c}{\(T_{\rm solve}\)} & \multicolumn{1}{c}{solve/PC} & \multicolumn{1}{c}{\(\eta_{\rm s}\)} & \multicolumn{1}{c}{memory} \\
\midrule
14.5 Hz & \(816\times816\times200\) & \revision{18} & 1 & \revision{0.005} s & \revision{3.462} s & \revision{0.192} s & 100.0\% & \revision{16.67} GiB \\
14.5 Hz & \(816\times816\times200\) & \revision{18} & 2 & \revision{0.002} s & \revision{1.847} s & \revision{0.103} s & \revision{93.7}\% & \revision{17.34} GiB \\
14.5 Hz & \(816\times816\times200\) & \revision{18} & 4 & \revision{0.002} s & \revision{1.148} s & \revision{0.064} s & \revision{75.4}\% & \revision{18.29} GiB \\
\addlinespace
29.1 Hz & \(1616\times1616\times384\) & \revision{36} & 8 & \revision{0.018} s & \revision{7.268} s & \revision{0.202} s & 100.0\% & \revision{126.10} GiB \\
29.1 Hz & \(1616\times1616\times384\) & \revision{36} & 16 & \revision{0.025} s & \revision{3.935} s & \revision{0.109} s & \revision{92.4}\% & \revision{130.62} GiB \\
29.1 Hz & \(1616\times1616\times384\) & \revision{36} & 32 & \revision{0.047} s & \revision{2.503} s & \revision{0.070} s & \revision{72.6}\% & \revision{138.39} GiB \\
\addlinespace
58.1 Hz & \(3216\times3216\times752\) & \revision{69} & 64 & \revision{0.076} s & \revision{14.566} s & \revision{0.211} s & 100.0\% & \revision{981.27} GiB \\
\bottomrule
\end{tabular}%
\end{adjustbox}
\end{table}

\revision{Table~}\revision{\ref{tab:overthrust-scaling}} \revision{confirms} \revision{that} the \revision{favorable} \revision{parallel}
\revision{performance} \revision{extends} to the geophysical model.  \revision{Strong-scaling} \revision{efficiency} \revision{is}
\revision{75.4}\revision{\%} \revision{from} \revision{one} \revision{to} \revision{four} \revision{GPUs} \revision{on} the \revision{native} \revision{mesh} and \revision{72.6}\revision{\%} \revision{from} \revision{eight} \revision{to}
\revision{thirty-two} \revision{GPUs} on \revision{the} \revision{next} \revision{finer} \revision{mesh,} \revision{consistent} with the \revision{synthetic}
\revision{benchmarks.}  \revision{Along} \revision{the} \revision{one-,} \revision{eight-,} and sixty-four-GPU \revision{refinement} \revision{sequence,}
the \revision{average} \revision{solve} time per preconditioner application \revision{stays} \revision{between} \revision{0.192}
\revision{and} \revision{0.211} \revision{seconds.}  \revision{The} \revision{finest} \revision{\(3216\times3216\times752\)} \revision{problem} \revision{is}
\revision{solved} \revision{in} \revision{14.6} \revision{seconds} \revision{at} \revision{\(58.1~\mathrm{Hz}\)}\revision{,} \revision{with} \revision{setup} \revision{taking} \revision{0.08}
\revision{seconds.}  \revision{These} \revision{results} \revision{demonstrate} \revision{that} \revision{the} \revision{independently} \revision{rediscretized}
\revision{physical} \revision{correction} \revision{and} \revision{shifted} \revision{auxiliary} \revision{inverse} \revision{remain} \revision{effective} \revision{when}
\revision{velocity} \revision{variation,} \revision{geological} \revision{interfaces,} \revision{and} \revision{increasing} \revision{frequency} \revision{are}
\revision{combined} \revision{in} \revision{one} \revision{large-scale} \revision{problem.}
\FloatBarrier

\revision{The CPU results extend this refinement sequence by one level: trilinear refinement to \(6401\times6401\times1473\) velocity samples (spacing \(3.125~\mathrm{m}\)) gives a \(6416\times6416\times1488\) mesh at \(116.2~\mathrm{Hz}\).  Table~\ref{tab:overthrust-weak-scaling} shows nearly constant solve/PC, with 97.1\% weak-scaling efficiency from 64 to 4096 cores, and Table~\ref{tab:overthrust-strong-scaling} shows that increasing the core count sixteenfold retains 66.2\% strong-scaling efficiency on the \(1616\times1616\times384\) mesh and 63.9\% on the \(3216\times3216\times752\) mesh, where the solve time drops to 43.4~s on 8192 cores; on the finest mesh, increasing the allocation from 4096 to 16384 and 32768 cores retains 86.1\% and 75.3\%, respectively, and reduces the solve time to 148.7~s.  The preconditioner applications are the same for every allocation of the first two meshes (36 and 69) and differ by at most one on the finest mesh (140 on 4096 cores, 139 otherwise).}

\begin{table}[H]
\centering
\caption{\revision{CPU} \revision{weak} \revision{scaling} for the \revision{SEG/EAGE} \revision{Overthrust} \revision{model.} \revision{Efficiency} \revision{is} \revision{normalized} \revision{by} \revision{the} \revision{actual} number of nodal unknowns \revision{and} \revision{measured} \revision{relative} \revision{to} the \revision{smallest} \revision{mesh} (8 MPI ranks of 8 OpenMP threads per node).}
\label{tab:overthrust-weak-scaling}
\footnotesize
\begin{adjustbox}{max width=\textwidth}
\begin{tabular}{ccccccccc}
\toprule
\multicolumn{1}{c}{\revision{frequency}} & \multicolumn{1}{c}{\revision{mesh}} & \multicolumn{1}{c}{\revision{CPU} \revision{cores}} & \multicolumn{1}{c}{\revision{PC} \revision{calls}} & \multicolumn{1}{c}{\revision{setup}} & \multicolumn{1}{c}{\revision{solve}} & \multicolumn{1}{c}{\revision{solve/PC}} & \multicolumn{1}{c}{\revision{\(\eta_{\rm w}\)}} & \multicolumn{1}{c}{\revision{peak/rank}} \\
\midrule
\revision{29.1 Hz} & \revision{\(1616\times1616\times384\)} & \revision{64} & \revision{36} & \revision{0.259} \revision{s} & \revision{236.746} \revision{s} & \revision{6.576} \revision{s} & \revision{100.0}\% & \revision{11.466} \revision{GiB} \\
\revision{58.1 Hz} & \revision{\(3216\times3216\times752\)} & \revision{512} & \revision{69} & \revision{0.265} \revision{s} & \revision{443.667} \revision{s} & \revision{6.430} \revision{s} & \revision{99.0}\% & \revision{11.179} \revision{GiB} \\
\revision{116.2 Hz} & \revision{\(6416\times6416\times1488\)} & \revision{4096} & \revision{140} & \revision{0.271} \revision{s} & \revision{902.276} \revision{s} & \revision{6.445} \revision{s} & \revision{97.1}\% & \revision{11.201} \revision{GiB} \\
\bottomrule
\end{tabular}%
\end{adjustbox}
\end{table}

\begin{table}[H]
\centering
\caption{\revision{CPU} \revision{strong} \revision{scaling} \revision{for} the \revision{SEG/EAGE} Overthrust \revision{model,} \revision{relative} to the \revision{smallest} \revision{allocation} \revision{for} \revision{each} \revision{mesh} (8 MPI ranks of 8 OpenMP threads per node).}
\label{tab:overthrust-strong-scaling}
\footnotesize
\begin{adjustbox}{max width=\textwidth}
\begin{tabular}{ccccccccc}
\toprule
\multicolumn{1}{c}{\revision{frequency}} & \multicolumn{1}{c}{\revision{mesh}} & \multicolumn{1}{c}{\revision{CPU} \revision{cores}} & \multicolumn{1}{c}{\revision{PC} \revision{calls}} & \multicolumn{1}{c}{\revision{setup}} & \multicolumn{1}{c}{\revision{solve}} & \multicolumn{1}{c}{\revision{solve/PC}} & \multicolumn{1}{c}{\revision{\(\eta_{\rm s}\)}} & \multicolumn{1}{c}{\revision{peak/rank}} \\
\midrule
\revision{29.1 Hz} & \revision{\(1616\times1616\times384\)} & \revision{64} & \revision{36} & \revision{0.259} \revision{s} & \revision{236.746} \revision{s} & \revision{6.576} \revision{s} & \revision{100.0}\% & \revision{11.466} \revision{GiB} \\
\revision{29.1 Hz} & \revision{\(1616\times1616\times384\)} & \revision{128} & \revision{36} & \revision{0.131} \revision{s} & \revision{121.621} \revision{s} & \revision{3.378} \revision{s} & \revision{97.3}\% & \revision{5.767} \revision{GiB} \\
\revision{29.1 Hz} & \revision{\(1616\times1616\times384\)} & \revision{256} & \revision{36} & \revision{0.067} \revision{s} & \revision{67.490} \revision{s} & \revision{1.875} \revision{s} & \revision{87.7}\% & \revision{3.006} \revision{GiB} \\
\revision{29.1 Hz} & \revision{\(1616\times1616\times384\)} & \revision{512} & \revision{36} & \revision{0.043} \revision{s} & \revision{37.105} \revision{s} & \revision{1.031} \revision{s} & \revision{79.8}\% & \revision{1.634} \revision{GiB} \\
\revision{29.1 Hz} & \revision{\(1616\times1616\times384\)} & \revision{1024} & \revision{36} & \revision{0.026} \revision{s} & \revision{22.336} \revision{s} & \revision{0.621} \revision{s} & \revision{66.2}\% & \revision{0.917} \revision{GiB} \\
\addlinespace
\revision{58.1 Hz} & \revision{\(3216\times3216\times752\)} & \revision{512} & \revision{69} & \revision{0.265} \revision{s} & \revision{443.667} \revision{s} & \revision{6.430} \revision{s} & \revision{100.0}\% & \revision{11.179} \revision{GiB} \\
\revision{58.1 Hz} & \revision{\(3216\times3216\times752\)} & \revision{1024} & \revision{69} & \revision{0.139} \revision{s} & \revision{227.715} \revision{s} & \revision{3.300} \revision{s} & \revision{97.4}\% & \revision{5.731} \revision{GiB} \\
\revision{58.1 Hz} & \revision{\(3216\times3216\times752\)} & \revision{2048} & \revision{69} & \revision{0.074} \revision{s} & \revision{126.981} \revision{s} & \revision{1.840} \revision{s} & \revision{87.3}\% & \revision{2.975} \revision{GiB} \\
\revision{58.1 Hz} & \revision{\(3216\times3216\times752\)} & \revision{4096} & \revision{69} & \revision{0.051} \revision{s} & \revision{71.697} \revision{s} & \revision{1.039} \revision{s} & \revision{77.4}\% & \revision{1.643} \revision{GiB} \\
\revision{58.1 Hz} & \revision{\(3216\times3216\times752\)} & \revision{8192} & \revision{69} & \revision{0.034} \revision{s} & \revision{43.424} \revision{s} & \revision{0.629} \revision{s} & \revision{63.9}\% & \revision{0.954} \revision{GiB} \\
\addlinespace
\revision{116.2 Hz} & \revision{\(6416\times6416\times1488\)} & \revision{4096} & \revision{140} & \revision{0.271} \revision{s} & \revision{902.276} \revision{s} & \revision{6.445} \revision{s} & \revision{100.0}\% & \revision{11.201} \revision{GiB} \\
\revision{116.2 Hz} & \revision{\(6416\times6416\times1488\)} & \revision{8192} & \revision{139} & \revision{0.146} \revision{s} & \revision{460.734} \revision{s} & \revision{3.315} \revision{s} & \revision{97.2}\% & \revision{5.747} \revision{GiB} \\
\revision{116.2 Hz} & \revision{\(6416\times6416\times1488\)} & \revision{16384} & \revision{139} & \revision{0.092} \revision{s} & \revision{260.149} \revision{s} & \revision{1.872} \revision{s} & \revision{86.1}\% & \revision{3.001} \revision{GiB} \\
\revision{116.2 Hz} & \revision{\(6416\times6416\times1488\)} & \revision{32768} & \revision{139} & \revision{0.047} \revision{s} & \revision{148.655} \revision{s} & \revision{1.070} \revision{s} & \revision{75.3}\% & \revision{1.572} \revision{GiB} \\
\bottomrule
\end{tabular}%
\end{adjustbox}
\end{table}

\FloatBarrier

\section{Conclusion}
\label{sec:conclusions}

\revision{We have developed a massively parallel three-grid preconditioner for
the three-dimensional high-frequency Helmholtz equation discretized by IOFD.
The method combines low-dispersion coarse-grid operators with a fully iterative
solution strategy, allowing accurate wave propagation to be resolved on compact
grids while retaining a matrix-free implementation on GPUs and CPUs.
Local Fourier analysis characterizes fine--coarse phase compatibility and
identifies admissible auxiliary parameters, providing guidance for the solver
configuration used throughout the numerical experiments.}

\revision{The numerical results demonstrate that this combination delivers
both wavefield accuracy and scalable solution of large problems.
Green-function comparisons confirm phase and amplitude accuracy over long
propagation distances at six points per shortest wavelength, while tests on
smooth, discontinuous, high-contrast, and geophysical media establish the
effectiveness of a common solver configuration.  The GPU and CPU experiments
further demonstrate strong and weak scaling, extending the calculations to
more than \revision{two trillion} unknowns and \revision{over 2200} wavelengths per coordinate
direction.  Together, these results show that dispersion-aware coarse
correction and a factorization-free hierarchy provide an effective approach
to large-scale high-frequency wave simulation.}

Future work will investigate block and recycled flexible Krylov methods for
multiple right-hand sides, together with alternative factorization-free
realizations of the \(2h\) inverse.  These extensions will preserve the
separation between dispersion control by the physical coarse operator and
stabilization of its inverse.

\section*{Declaration of AI Assistance}
The authors acknowledge the assistance of ChatGPT (OpenAI) in the theoretical
analysis, code implementation, and language editing of this manuscript.
The authors take full responsibility for the accuracy and integrity of all
content in the paper.

\bibliographystyle{siamplain}
\bibliography{olfd_references}

\end{document}